\documentclass[]{siamltex}
\usepackage{amsmath}
\usepackage{amssymb}
\usepackage{graphicx}

\usepackage{tikz}
\usepackage{pgfplots}
\usepackage[ComplainAboutUnexternalized]{externaltikz}
\usepackage{calc}
\usepackage{subcaption}

\usepackage{hyperref}
 \usepackage{mathtools}
 \mathtoolsset{showonlyrefs}

\newcommand{\uuK}{\uu_K} 
\newcommand{\uuB}{\uu_B} 
\newcommand{\uu}{u} 
\newcommand{\uuhat}{\hat{u}} 

\newtheorem{remark}{Remark}

\begin{document}

\title{Kinetic layers and coupling conditions for  nonlinear scalar  equations on networks}

\author{R. Borsche\footnotemark[1] 
     \and  A. Klar\footnotemark[1] \footnotemark[2]}
\footnotetext[1]{Technische Universit\"at Kaiserslautern, Department of Mathematics, Erwin-Schr\"odinger-Stra{\ss}e, 67663 Kaiserslautern, Germany (\{borsche, klar\}@mathematik.uni-kl.de)}
\footnotetext[2]{Fraunhofer ITWM, Fraunhoferplatz 1, 67663 Kaiserslautern, Germany} 
 
\date{}


\maketitle

\begin{abstract}
We consider a kinetic relaxation  model and an associated macroscopic scalar nonlinear hyperbolic equation on a network. 
Coupling conditions for the macroscopic equations are derived from the kinetic coupling conditions via an asymptotic analysis near the nodes of the network. 
This analysis leads to the combination of kinetic half-space problems with Riemann problems at the junction.
Detailed numerical comparisons between the different models show the agreement of the coupling conditions for the case of tripod networks.
\end{abstract}




\section{Introduction}

Coupling conditions for macroscopic partial differential equations on networks including, for example,  diffusion equations,  wave type equations or Euler equations have been discussed in many papers, see, for example,  \cite{Cor17,BNR14,GPBook,LS02,CGP05,BHK06a,CG08,CM08}.
Coupling conditions for the underlying kinetic equations on networks have been considered in a much smaller number of publications \cite{HM09,BKP16}. 
In \cite{BKP16} a first attempt to derive a coupling condition for a
macroscopic equation from the underlying kinetic model has been presented for the case of the 
 kinetic  chemotaxis equation. A more   general and more accurate procedure to derive coupling conditions for macroscopic 
equations from the underlying kinetic ones has been  discussed for linear systems in \cite{BK17} using  an asymptotic analysis of the situation near the nodes.
It has been  motivated by the classical procedure to find kinetic slip boundary conditions for macroscopic equations derived from underlying kinetic equations based on the analysis of the 
kinetic layer, see \cite{BSS84,BLP79,G08,UTY03} for kinetic equation or 
\cite{WY99,WX99,LX96,X04} for the case of hyperbolic relaxation systems. At each node of the network a  fixpoint problem involving the coupled solution of  half-space problems for each edge  has been approximated. 
Using such a procedure explicit coupling conditions have been derived for  the linear wave equation from an underlying linear kinetic model in \cite{BK17}.
 
In the present  work, we  extend this analysis and derive coupling conditions for nonlinear scalar  equations 
on a network from an underlying  kinetic relaxation model. We concentrate on  a two equation kinetic relaxation model leading in the limit to the Burger's equation.

The paper is organized in the following way.  In section \ref{equations} we present the relaxation model and the scalar conservation law. 
In section \ref{layeranalysis}  kinetic boundary layers are discussed, as well as the combination of these layer solutions with suitable Riemann solvers. This leads to classical boundary conditions for the Burgers problem depending on the kinetic boundary condition.
 
In the following section \ref{couplingconditions}  coupling conditions for the scalar hyperbolic problem are discussed
and derived from the kinetic coupling conditions. We start  with a  node with 2 edge.
In this case the above procedure leads simply to the solution of a Riemann problem at the junction without a kinetic layer.
This changes in the case of a  node with three edges. For this case we  derive in section \ref{3node} explicit coupling conditions for the macroscopic equation
based on the kinetic coupling conditions.

Finally, the solution of the macroscopic equations on the network are numerically compared  to the full solutions of the kinetic equation on the network  in section \ref{Numerical results}.

\section{Equations}
\label{equations}

We consider  the following   relaxation  model for   $x \in \mathbb{R} , v_1<  0 < v_2 $ and $F= F(\uu)$.
\begin{align}\label{bgk}
\begin{aligned}
\partial_t f_1 + v_1 \partial_x f_1 = -\frac{1}{\epsilon} \left(f_1-\frac{v_2 \uu-F(\uu)}{v_2-v_1} \right) \\
\partial_t f_2 + v_2 \partial_x f_2 = -\frac{1}{\epsilon} \left(f_2-\frac{ F(\uu)-v_1 \uu}{v_2-v_1}\right)
\end{aligned}
\end{align}
with
$
\uu = f_1+f_2
$. 
Defining the flux
$
\uuhat=v_1 f_1 + v_2 f_2
$
yields
\begin{align*}
f_1 = \frac{v_2 \uu-\uuhat}{v_2 -v_1}
\ ,\quad 
f_2= \frac{\uuhat-v_1 \uu}{v_2 -v_1}\ .
\end{align*}
The associated macroscopic equation  for $\epsilon \rightarrow 0$ is 
a conservation law for the quantity  $u$ given by 
\begin{align}\label{eq:CL}
\begin{aligned}
\partial_t \uu + \partial_x F(\uu) =0\ .
\end{aligned}
\end{align}

Rewriting  the kinetic problem in terms of  $\uu $ and $\uuhat$ yields
\begin{align}\label{euler}
\begin{aligned}
\partial_t \uu + \partial_x \uuhat &=0\\
\partial_t \uuhat + \partial_x  P(\uu,\uuhat ) &=-\frac{1}{\epsilon} \left(\uuhat-F(\uu) \right) \ .
\end{aligned}
\end{align}
with
\begin{align}\label{eq:P}
P(\uu,\uuhat) = (v_1+v_2) \uuhat-v_1 v_2 \uu\ .
\end{align}

Convergence of the kinetic equation as $\epsilon \rightarrow 0$ is obtained under the subcharacteristic condition \cite{Liu}
\begin{align}\label{eq:subchar}
v_1 \le F^\prime (\uu) \le v_2\ .
\end{align}

The boundary value problem for \eqref{bgk} is simple, since the equations are linear on the left hand side with $v_1<0<v_2$.
Thus, we have to prescribe $f_2$ at the left boundary and $f_1$ at the right boundary.
For the nonlinear hyperbolic limit problem, boundary conditions have been considered in many works \cite{BRN,O96,AS15,ACD17,BFK07}.

For boundary conditions and  layers of hyperbolic problems with stiff relaxation terms and for the derivation of conditions for the corresponding limit equations we refer for example to \cite{WY99,WX99,LX96,X04}.

We determine the boundary conditions for the limit equations via a combination 
of an analysis of the kinetic layer with the solution of a Half-Riemann problem for the limit equation, compare  for example \cite{WX99,AM04}.
A similar procedure will then be used to find kinetic based  coupling conditions for the Burgers equation on a network.

To proceed, we first state the layer equations.
We consider the left boundary of the domain located at $x=0$.
A rescaling of the spacial coordinate near the boundary with $x \rightarrow \frac{x}{\epsilon }$ gives the layer problem on $[0,\infty)$
as
\begin{align}\label{layer}
\begin{aligned}
 v_1 \partial_x f_1 &= \frac{v_2 \uu-F(\uu)}{v_2-v_1} -f_1 \\
 v_2 \partial_x f_2 &= \frac{ F(\uu)-v_1 \uu}{v_2-v_1} -f_2 \ .
\end{aligned}
\end{align}

In the macroscopic variables $\uu,\uuhat$ this is 
\begin{align*}
\begin{aligned}
 \partial_x \uuhat &=0\\
\partial_x  P(\uu,\uuhat ) &=F(\uu) -\uuhat\ ,
\end{aligned}
\end{align*}
with $P$ as in \eqref{eq:P}.
This gives
\begin{align*}
\begin{aligned}
 \partial_x \uuhat &=0\\
-v_1 v_2 \partial_x  \uu &=F(\uu) -\uuhat\ 
\end{aligned}
\end{align*}
and therefore
\begin{align}\label{layereq}
\begin{aligned}
\uuhat &= C = const\\
a  \partial_x \uu &= F(\uu) -C 
\end{aligned}
\end{align}
with 
$a= - v_1 v_2 >0$.

\begin{remark}
For a right boundary  we obtain the layer problem as
\begin{align*}
- a  \partial_x \uu = F(\uu) -C\ .
\end{align*}
\end{remark}

For the following computations we concentrate on 
 the Burgers equations and choose  $F(\uu) = \uu^2$.
 Thus, the corresponding subcharacteristic condition \eqref{eq:subchar} is $v_1 \le 2\min (u)$, $v_2 \ge 2 \max (u)$.

\section{Boundary conditions}
\label{layeranalysis}

In a first step we determine the boundary conditions for the hyperbolic limit equation from the kinetic boundary condition. 
We use the boundary layer equations and couple them with Half-Riemann solvers. 
First we discuss the solution in the boundary layer.

\subsection{The boundary layer equation}
The boundary layer equation  near a left boundary is given by 
$$
 a \partial_x \uu = \uu^2 -C.
$$
For  $C>0$ this problem has  two fixpoints $\uu=\pm \sqrt{C}$, where 
$\sqrt{C}$ is instable and  $-\sqrt{C}$ is a stable fixpoint.
The domain of attraction of the stable fix-point is $(-\infty,\sqrt{C})$.

The explicit solution is given by
\begin{align*}
\uu(x) =  \sqrt{C} \mbox{tanh} (-  \sqrt{C} (x+C_2)/a) 
&& \text{ for }\quad \vert \uu(0)\vert < \sqrt{C}
\ 
\end{align*}
and 
\begin{align*}
\uu(x) =\sqrt{C} \mbox{coth} (-  \sqrt{C} (x+C_2)/a)  
&& \text{ for }\quad \vert \uu(0)\vert > \sqrt{C}
\ .
\end{align*}
We determine $C_2$ from
\begin{align*}
\uu(0) = \sqrt{C} \mbox{tanh} (-  \sqrt{C} C_2/a) 
&& \text{ for }\quad \vert \uu(0)\vert < \sqrt{C}\ 
\end{align*}
and
\begin{align*}
\uu(0) =\sqrt{C} \mbox{coth} (-  \sqrt{C} C_2/a) 
&& \text{ for }\quad \vert \uu(0)\vert > \sqrt{C}\ .
\end{align*}
One observes that for $\uu(0) < \sqrt{C}$ the limit $x \rightarrow \infty$ leads to $\uu(x ) \rightarrow  - \sqrt{C}$
and for $\uu(0) > \sqrt{C}$ the layer solution diverges at $x=-C_2 = \frac{a}{\sqrt{C}} \mbox{arcoth}(\frac{\uu(0)}{\sqrt{C}})$.

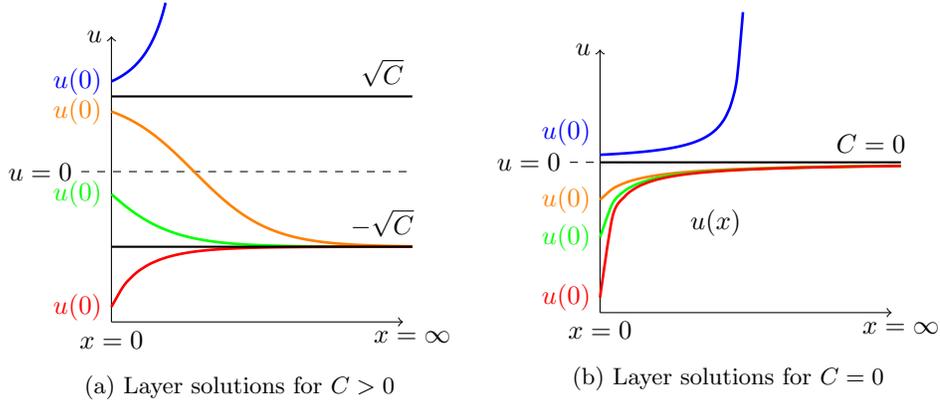
\begin{figure}[h]
	\externaltikz{Layer_solutions}{
	\begin{subfigure}{0.49\textwidth}
	\begin{tikzpicture}
	\def\len{4}
	\def\low{-2}
	
	\def \uuzero {0.8}
	\def \uuzeroo {-0.3}
	\def \uuzerotw {-1.8}
	\def \uuzerotr {1.2}
	\def \Ckonst {1}
	\def \sCkonst {1}
	\def \Ctwo {-1/(2*\sCkonst)*ln((\sCkonst+\uuzero)/(\sCkonst-\uuzero))}
	\def \Ctwoo {-1/(2*\sCkonst)*ln((\sCkonst+\uuzeroo)/(\sCkonst-\uuzeroo))}
	\def \Ctwotw {-1/(2*\sCkonst)*ln((\uuzerotw+\sCkonst)/(\uuzerotw-\sCkonst))}
	\def \Ctwotr {-1/(2*\sCkonst)*ln((\uuzerotr+\sCkonst)/(\uuzerotr-\sCkonst))}
	
	\node[] (xinf) at (\len,\low){};
	\node[below] at (0,\low){$x=0$};
	\node[below] at (\len,\low){$x=\infty$};
	\draw[dashed] (-0.1*\len,0)--(\len,0) node[left,pos = 0]{$\uu=0$};
	\draw[->] (0,\low)--(xinf);
	\draw[->] (0,\low)--(0,1.8) node[left]{$\uu$};
	\draw[orange,line width=1pt,domain=0.0:\len,smooth,variable=\x,] plot ({\x},{\sCkonst*tanh(-\sCkonst*(\x+\Ctwo))}) ;
	\draw[green,line width=1pt,domain=0.0:\len,smooth,variable=\x,] plot ({\x},{\sCkonst*tanh(-\sCkonst*(\x+\Ctwoo))}) ;
	\draw[red,line width=1pt,domain=0.0:\len,smooth,variable=\x,] plot ({\x},{\sCkonst*cosh(-\sCkonst*(\x+\Ctwotw))/sinh(-\sCkonst*(\x+\Ctwotw))}) ;
	\draw[blue,line width=1pt,domain=0.0:0.18*\len,smooth,variable=\x,] plot ({\x},{\sCkonst*cosh(-\sCkonst*(\x+\Ctwotr))/sinh(-\sCkonst*(\x+\Ctwotr))}) ;
	\draw[thick] (0,\sCkonst)--(\len,\sCkonst) node[above,pos = 0.9]{$\sqrt{C}$};
	\draw[thick] (0,-\sCkonst)--(\len,-\sCkonst) node[above,pos = 0.9]{$-\sqrt{C}$};
	\node[orange,left] at (0.*\len,\uuzero) {$\uu(0)$};
	\node[green,left] at (0.*\len,\uuzeroo) {$\uu(0)$};
	\node[red,left] at (0.*\len,\uuzerotw) {$\uu(0)$};
	\node[blue,left] at (0.*\len,\uuzerotr) {$\uu(0)$};
	\end{tikzpicture}
	\caption{Layer solutions for $C >0$}
	\end{subfigure}
	\begin{subfigure}{0.49\textwidth}
	\begin{tikzpicture}
	\def\len{4}
	\def\low{-2}
	
	\def \uuzero {-0.5}
	\def \uuzeroo {-1}
	\def \uuzerotw {-1.8}
	\def \uuzerotr {0.1}
	
	\node[] (xinf) at (\len,\low){};
	\node[below] at (0,\low){$x=0$};
	\node[below] at (\len,\low){$x=\infty$};
	\draw[dashed] (-0.1*\len,0)--(\len,0) node[left,pos = 0]{$\uu=0$};
	\draw[->] (0,\low)--(xinf);
	\draw[->] (0,\low)--(0,1.5) node[left]{$\uu$};
	\draw[orange,line width=1pt,domain=0.0:\len,smooth,variable=\x,] plot ({\x},{-1/(\x/\len*20-1/(\uuzero))}) ;
	\draw[green,line width=1pt,domain=0.0:\len,smooth,variable=\x,] plot ({\x},{-1/(\x/\len*20-1/(\uuzeroo))}) ;
	\draw[red,line width=1pt,domain=0.0:\len,smooth,variable=\x,] plot ({\x},{-1/(\x/\len*20-1/(\uuzerotw))}) ;
	\draw[blue,line width=1pt,domain=0.0:1.9,smooth,variable=\x,] plot ({\x},{-1/(\x/\len*20-1/(\uuzerotr))}) ;
	\draw[thick] (0,0)--(\len,0) node[above,pos = 0.9]{$C=0$};
	\node[orange,left] at (0,\uuzero) {$\uu(0)$};
	\node[green,left] at (0,\uuzeroo) {$\uu(0)$};
	\node[red,left] at (0,\uuzerotw) {$\uu(0)$};
	\node[blue,left] at (0,\uuzerotr+0.3) {$\uu(0)$};
	\node[left] at (0.5*\len,-0.8) {$\uu(x)$};
	\end{tikzpicture}
	\caption{Layer solutions for $C =0$}
	\end{subfigure}
	}
	\caption{Possible solutions to the layer equation.}
	\label{fig:SketchLayer}
\end{figure}

For $C=0$ we obtain
\begin{align*}
\uu(x) = - \frac{a}{x+C_2} =- \frac{a}{x-\frac{a}{\uu(0)}} 
\end{align*}
and convergence to $0$ for $\uu(0) <0$ and divergence for $\uu(0) >0$.
The solutions are sketched in figure \ref{fig:SketchLayer}.

In the following we use the the notation $(U)$ for the unstable solution $\uu(x) = \sqrt{C}$ and  the notation $(S)$
for the (partially) stable solutions. The asymptotic states as $x \rightarrow \infty $ are denoted by $u_K$. 
The layer solution for the right boundary can be discussed analogously.

\subsection{Riemann Problem}
Since the layer solution can not cover the full range of possible states at a boundary, we have to consider additionally a Riemann Problem for the Burgers equation connecting the state in the domain with the layer.
In particular, for the left boundary we need to know, which asymptotic states $\uuK$ from the kinetic layer  can be connected to a given right side state from the Burgers equation $\uu_B$ using only waves with non-negative speeds.
For the Burgers equation  we have the following cases:
\begin{enumerate}
	\item[{\bf RP1}]
	$\uu_B\geq 0$ $\Rightarrow$ $\uuK\in [0,\infty)$, since there is either an arbitrary wave with positive speed, if $\uuK>0$, or a rarefaction wave starting at $\uu=0$.
	\item[{\bf RP2}]
	$\uu_B<0$ $\Rightarrow$ $\uuK\in \{\uu_B\}\cup (-\uu_B,\infty)$, since
	there is either 
	no wave or a  shock wave moving to the right.
\end{enumerate}
Thus, for a given $\uu_B$ we can select $\uuK$ only from the above subsets.
For a boundary on the right hand, we study the analogous cases.

\subsection{Macroscopic boundary conditions}\label{sec:MacroBoundaryCond}
To find the macroscopic boundary conditions at a left boundary from the underlying kinetic problem, we combine  the solution of the half space problem \eqref{layereq} on $[0,\infty]$ with asymptotic solution $\uu_K$ and of a Riemann Problem  with left state $\uu_K$  and right state $\uu_B$.

We consider  a domain  $x \ge 0$ and determine macroscopic boundary conditions at $x=0$ in the following way.
For the kinetic problem we prescribe $f_2(x=0)$. Moreover, the actual macroscopic value at $x=0^+$ is denoted by $\uu_B$. From these two values we have to determine a (potentially new) boundary value for the macroscopic solution $\uu_K$ and a value for  $f_1(0)$,
the outgoing kinetic value.
We consider different cases coupling stable or unstable layer solutions, denoted by $(U)$ or $(S)$, and Riemann problem solutions $RP1$ or $RP2$.

\paragraph{{\bf Case 1, RP1-U}} $\uu_B >0$ and $f_2(0) >0$. 

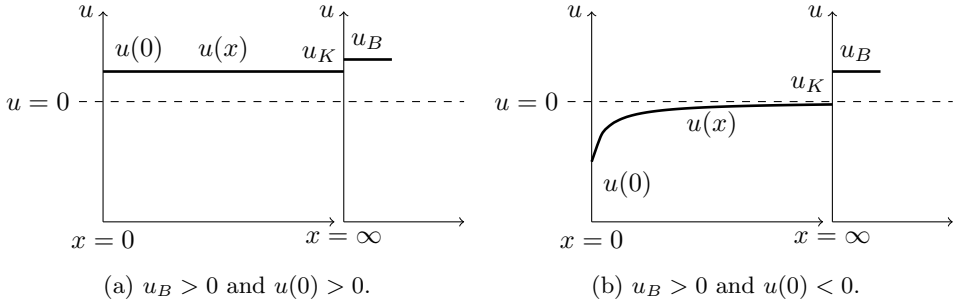
\begin{figure}[h]
	\externaltikz{BLayer_uBpos}{
	\begin{subfigure}{0.49\textwidth}	
	\begin{tikzpicture}[scale = 0.8]
	\def\len{4}
	\def\low{-2}
	
	\def \uuzero {-1}
	\def \uuB {0.7}
	\def \uuK {0.5}
	\def \uuk {-0.75}
	\def \Ckonst {\uuk*\uuk}
	\def \Ctwo {-1/(2*\uuk)*ln((\uuk+\uuzero)/(\uuk-\uuzero))}
	
	\node[] (xinf) at (\len,\low){};
	\node[below] at (0,\low){$x=0$};
	\node[below] at (\len,\low){$x=\infty$};
	\draw[dashed] (-0.1*\len,0)--(1.5*\len,0) node[left,pos = 0]{$\uu=0$};
	\draw[->] (0,\low)--(xinf);
	\draw[->] (\len,\low)--(1.5*\len,\low);
	\draw[->] (0,\low)--(0,1.5) node[left]{$\uu$};
	\draw[->] (\len,\low)--(\len,1.5) node[left]{$\uu$};
	\draw[line width=1pt] (\len,\uuB)--(1.2*\len,\uuB) node[above,pos=0.5] {$\uu_B$};
	\draw[line width=1pt,domain=0.0:\len,smooth,variable=\x,] plot ({\x},{\uuK}) ;
	\node[above] at (0.15*\len,\uuK) {$\uu(0)$};
	\node[above] at (0.5*\len,\uuK) {$\uu(x)$};
	\node[above] at (\len*0.9,\uuK) {$\uu_K$};	
	\end{tikzpicture}
	\caption{
		 $\uu_B>0$ and $\uu(0)>0$.
		}
	\end{subfigure}
	\begin{subfigure}{0.49\textwidth}
	\begin{tikzpicture}[scale = 0.8]
	\def\len{4}
	\def\low{-2}
	
	\def \uuzero {-1}
	\def \uuB {0.5}
	\def \uuk {0}
	\def \Ckonst {\uuk*\uuk}
	\def \Ctwo {-1/(2*\uuk)*ln((\uuk+\uuzero)/(\uuk-\uuzero))}

	\node[] (xinf) at (\len,\low){};
	\node[below] at (0,\low){$x=0$};
	\node[below] at (\len,\low){$x=\infty$};
	\draw[dashed] (-0.1*\len,0)--(1.5*\len,0) node[left,pos = 0]{$\uu=0$};
	\draw[->] (0,\low)--(xinf);
	\draw[->] (\len,\low)--(1.5*\len,\low);
	\draw[->] (0,\low)--(0,1.5) node[left]{$\uu$};
	\draw[->] (\len,\low)--(\len,1.5) node[left]{$\uu$};
	\draw[line width=1pt] (\len,\uuB)--(1.2*\len,\uuB) node[above,pos=0.5] {$\uu_B$};
	\draw[line width=1pt,domain=0.0:\len,smooth,variable=\x,] plot ({\x},{-1/(\x/\len*20-1/(\uuzero))}) ;
	\node[below] at (0.15*\len,\uuzero) {$\uu(0)$};
	\node[below] at (0.5*\len,\uuk) {$\uu(x)$};
	\node[above] at (\len*0.9,\uuk) {$\uu_K$};	
	\end{tikzpicture}
	\caption{$\uu_B>0$ and $\uu(0)<0$. }
	\end{subfigure}
	}
	
	\caption{Boundary layer and Riemann problem solution for positive $\uu_B$.}
	\label{positive}	
\end{figure}

The layer solution is
$$
\uu(x) =  \sqrt{C} >0\ .
$$
Determine  the value of $C > 0$ from  
$$
f_2 (0) = \frac{\uuhat-v_1 \uu(0)}{v_2-v_1} = \frac{C-v_1 \sqrt{C}}{v_2-v_1}\ .
$$ 
This equation has a positive solution $C$ under the above condition on $f_2(0)$.
This gives 
$$
C =  \frac{1}{4}\left( v_1 + \sqrt{v_1^2+4 (v_2-v_1)  f_2(0)} \right)^2
$$ 
and the new boundary condition $\uuK=\sqrt{C}$ given by the value of the layer solution at $\infty$.
The outgoing layer solution is $f_1 (0) = \frac{C- v_1 \sqrt{C}}{v_2-v_1} $.
This yields
$$f_1 (0) = \sqrt{C} - f_2(0)$$
and
$\uu(0) = \sqrt{C}$. See Figure \ref{positive}-(a).

\paragraph{{\bf Case 2, RP1-S}}  $\uu_B >0$ and $f_2(0) <0$.

In this case  $C =0$ and 
$\uu(0)$ is given by
$
f_2(0) = \frac{C-v_1 \uu(0)}{v_2-v_1}
$
or
$$
\uu(0) = \frac{ -f_2(0) (v_2-v_1)}{v_1} <0
$$
The layer solution is then
$
\uu(x) =- \frac{a}{x- a/\uu(0)}
$.
We do not need the exact form of the solution, but only the fact that for all $\uu(0) <0$ the asymptotic value is $0$.
The new boundary condition is $\uuK =0$
and 
$$
f_1(0) = -f_2(0)\frac{v_2}{v_1}.
$$
In this case $\uu(0)$ and $\uu_B$ cannot be connected by an outgoing Burgers wave.
The kinetic layer solution takes care for a  part (from $\uu(0)$ to $0$) of the full jump from $\uu(0)$ to $\uu_B$, see Figure \ref{positive}-(b).

\paragraph{{\bf Case 3, RP2-S}} $\uu_B <0$ and $f_2(0) \le \frac{\uu_B^2 +v_1\uu_B}{v_2-v_1}$.

In this case 
the value of the layer solution at infinity is given by $\uu_B$. Thus,
$$
\uu(\infty) =- \sqrt{C} = \uu_B\ .
$$
Therefore  $C = \uu_B^2$.
Determine  $\uu(0) $ from $
f_2 (0)= \frac{C-v_1 \uu(0)}{v_2-v_1} = \frac{\uu_B^2-v_1 \uu(0)}{v_2-v_1}
$,
i.e. under the above assumption on $f_2(0)$ we have 
$$
\uu(0) = \frac{\uu_B^2}{v_1} - f_2(0) \frac{v_2-v_1}{v_1} \le -\uu_B = \sqrt{C}\ .
$$

Then, the  layer solution is given by the formulas in the last subsection and converges to the stable fixpoint $-\sqrt{C}$. Moreover,
$$
f_1(0) = \uu(0) -f_2(0) = \frac{\uu_B^2}{v_1} - f_2(0) \frac{v_2}{v_1}\ .
$$

\begin{figure}[h]
	\externaltikz{BLayer_uuBneg}{
	\begin{subfigure}{0.49\textwidth}
	\begin{tikzpicture}[scale = 0.8]
	\def\len{4}
	\def\low{-2}
	
	\def \uuzero {0.4}
	\def \uuB {-0.5}
	\def \uuk {-0.5}
	\def \Ckonst {\uuk*\uuk}
	\def \Ctwo {-1/(2*\uuk)*ln((\uuk+\uuzero)/(\uuk-\uuzero))}
	
	\node[] (xinf) at (\len,\low){};
	\node[below] at (0,\low){$x=0$};
	\node[below] at (\len,\low){$x=\infty$};
	\draw[dashed] (-0.1*\len,0)--(1.5*\len,0) node[left,pos = 0]{$\uu=0$};
	\draw[->] (0,\low)--(xinf);
	\draw[->] (\len,\low)--(1.5*\len,\low);
	\draw[->] (0,\low)--(0,1.5) node[left]{$\uu$};
	\draw[->] (\len,\low)--(\len,1.5) node[left]{$\uu$};
	\draw[line width=1pt] (\len,\uuB)--(1.2*\len,\uuB) node[above,pos=0.5] {$\uu_B$};
	\draw[line width=1pt,domain=0.0:\len,smooth,variable=\x,] plot ({\x},{\uuk*tanh(-\uuk*(\x/\len*10+\Ctwo))});
	\node[above] at (0.15*\len,{\uuk*tanh(-\uuk*\Ctwo)}) {$\uu(0)$};
	\node[below] at (0.5*\len,\uuB) {$\uu(x)$};
	\node[above] at (\len*0.9,\uuB) {$\uu_K$};	
	\end{tikzpicture}
	\caption{$\uu_B<0$ and   $\uu(0) \le -\uu_B$. }
	\end{subfigure}
	\begin{subfigure}{0.49\textwidth}
	\begin{tikzpicture}[scale = 0.8]
	\def\len{4}
	\def\low{-2}
	
	\def \uuzero {-1}
	\def \uuB {-0.5}
	\def \uuk {0.75}
	\def \Ckonst {\uuk*\uuk}
	\def \Ctwo {-1/(2*\uuk)*ln((\uuk+\uuzero)/(\uuk-\uuzero))}
	
	\node[] (xinf) at (\len,\low){};
	\node[below] at (0,\low){$x=0$};
	\node[below] at (\len,\low){$x=\infty$};
	\draw[dashed] (-0.1*\len,0)--(1.5*\len,0) node[left,pos = 0]{$\uu=0$};
	\draw[->] (0,\low)--(xinf);
	\draw[->] (\len,\low)--(1.5*\len,\low);
	\draw[->] (0,\low)--(0,1.5) node[left]{$\uu$};
	\draw[->] (\len,\low)--(\len,1.5) node[left]{$\uu$};
	\draw[line width=1pt] (\len,\uuB)--(1.2*\len,\uuB) node[above,pos=0.5] {$\uu_B$};
	\draw[line width=1pt,domain=0.0:\len,smooth,variable=\x,] plot ({\x},{\uuk}) ;
	\node[above] at (0.15*\len,\uuk) {$\uu(0)$};
	\node[above] at (0.5*\len,\uuk) {$\uu(x)$};
	\node[above] at (\len*0.9,\uuk) {$\uu_K$};	
	\end{tikzpicture}
	\caption{$\uu_B<0$ and   $\uu(0) \ge -\uu_B$. }
	\end{subfigure}
	}
	\caption{Boundary layer and Riemann problem solution for negative $\uu_B$.}
	\label{negative}
\end{figure}
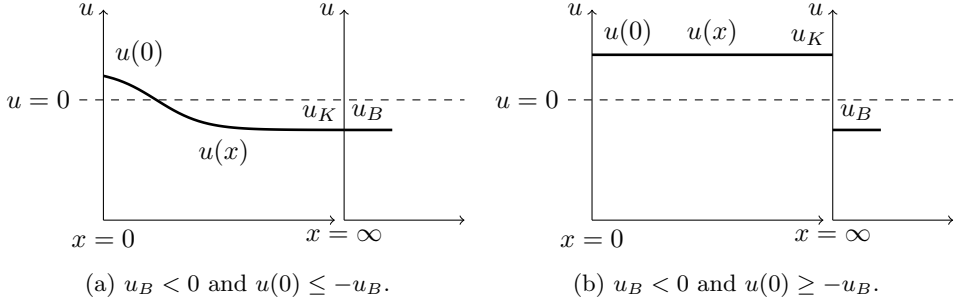

In this case $\uu(0)$ and $\uu_B$ cannot be connected by an outgoing Burgers wave.
The kinetic layer solution handles the full jump, see Figure \ref{negative}-(a).

\paragraph{{\bf Case 4, RP2-U}} $\uu_B <0$ and $f_2(0) \ge \frac{\uu_B^2 +v_1\uu_B}{v_2-v_1}$.
As in the first case  the layer solution is given by $$
\uu(x) =  \sqrt{C} \ .
$$
The value of $C > 0$ is determined from  $$
f_2 (0) = \frac{C-v_1 \sqrt{C}}{v_2-v_1}
$$ 
or
$$
C =  \frac{1}{4}\left( v_1 + \sqrt{v_1^2+4 (v_2-v_1)  f_2(0)} \right)^2\ .
$$ 
The boundary condition $\uuK=\sqrt{C}$ and the outgoing layer solution are given as in case 1.
In this case the full jump can be covered by a Riemann problem solution, see Figure \ref{negative}-(b).

\subsection{Summary boundary conditions}
We conclude this section with a short summary of the boundary conditions at the left and right boundaries. We simplify the situation assuming $v_2 = -v_1=v$.

\subsubsection{ Left boundary}

\paragraph{{\bf Case 1, RP1-U and RP2-U}} (ingoing flow) $\uu_B >0$ and $f_2(0)>0$ or  $\uu_B <0$ and $f_2(0) \ge \frac{\frac{\uu_B^2}{v} -\uu_B}{2}$
Then
\begin{align*}
\uuK=\frac{v}{2} \left(-1 + \sqrt{1+8  \frac{f_2(0)}{v}}\right)
\ , &&f_1(0) = \uuK - f_2(0)\ .
\end{align*}

\paragraph{{\bf Case 2, RP1-S}} (transsonic flow) $\uu_B >0$ and $f_2(0)<0$.
Then
\begin{align*}
\uuK=0\ ,&&f_1(0) = f_2(0)\ .
\end{align*}

\paragraph{{\bf Case 3, RP2-S}}  (outgoing flow) $\uu_B <0$ and $f_2(0) \le \frac{\frac{\uu_B^2}{v} -\uu_B}{2}$
Then
\begin{align*}
\uuK= \uu_B
\ ,&&
f_1(0) =  f_2(0) -\frac{\uu_B^2}{v} \ .
\end{align*}

\subsubsection{ Right boundary} 

\paragraph{{\bf Case 1, RP1-U and RP2-U}} (ingoing flow) $\uu_B <0$ and $f_1(0)<0$ or  $\uu_B >0$ and $f_1(0) \le -\frac{\uu_B^2}{2v} -\frac{\uu_B}{2}$.
Then
\begin{align*}
\uuK=\frac{v}{2} \left(1 - \sqrt{1-8  \frac{f_1(0)}{v}}\right)
\ ,&&
f_2(0) = \uuK - f_1(0)\ .
\end{align*}

\paragraph{{\bf Case 2, RP1-S}}  (transsonic flow) $\uu_B <0$ and $f_1(0)>0$.
Then
\begin{align*}
	\uuK=0\ ,&&
	f_2(0) = f_1(0)\ .	
\end{align*}

\paragraph{{\bf Case 3, RP2-S}}  (outgoing flow) $\uu_B >0$ and $f_1(0) \ge -\frac{\uu_B^2}{2v} -\frac{\uu_B}{2}$
Then,
\begin{align*}
	\uuK=\uu_B\ ,&&
	f_2(0) =  f_1(0) +\frac{\uu_B^2}{v} \ .
\end{align*}

\begin{remark}
Thus, we obtain conditions similar to the ones obtained in \cite{AM04}. For a proof of convergence of the kinetic BVP to the macroscopic one, we refer  \cite{WX99}, where the transonic case has been excluded. 
\end{remark}

\section{Coupling conditions}
\label{couplingconditions}

In this and the following section we consider the kinetic and Burgers problems on a network. In particular, we consider nodes with two and three edges.
In the present case the orientation of the edges is important. 
In the kinetic case we have to prescribe at the node for each ingoing edge the quantity $f_1$ and for the outgoing edges $f_2$. 
The macroscopic coupling conditions are then derived from the kinetic conditions using the above procedures.
We start the investigation by collecting the different cases and considering a simple node with only one incoming and one outgoing edge in this section and the case of a node with 3 edges in the following section.

\subsection{Summary of possible cases}

Before considering  the coupling problem for the different nodes in the next sections  we  collect again all combinations of possible  cases for the layer and half-Riemann problems.

\subsubsection{Half-Riemann problems}

\paragraph{The half-Riemann Problem at the left boundary}

		\begin{align*}
		\uu_B&\geq 0 \ (\text{RP 1})  \quad &\Rightarrow\quad \uuK &\in [0,\infty)\\
		\uu_B&< 0\ (\text{RP 2})  \quad &\Rightarrow\quad \uuK &\in \{\uu_B\}\cup [-\uu_B,\infty)\\
		\end{align*}		
\paragraph{The half-Riemann Problem at the right boundary}		
		\begin{align*}
		\uu_B&\leq 0\ (\text{RP 1})  \quad &\Rightarrow\quad \uuK &\in (-\infty,0]\\
		\uu_B&> 0\ (\text{RP 2})  \quad &\Rightarrow\quad \uuK &\in (-\infty,-\uu_B]\cup\{\uu_B\}\\
		\end{align*}

\subsubsection{Layer  problems}

\paragraph{The Layer Problem at the left boundary} 

			\begin{align*}
			&\quad \uuK = \sqrt{C} \quad \Rightarrow\quad \uu (0) = \sqrt{C},\ \uuhat = C>0
			&\quad \text{(U)}\\
			&\left.\begin{array}{ll}
				\uuK = 0 \quad &\Rightarrow\quad  \uu(0) \in (-\infty,0],\ \uuhat = 0\\
				\uuK = -\sqrt{C} \quad &\Rightarrow\quad \uu(0) \in (-\infty,\sqrt{C}),\ \uuhat = C>0
			\end{array}\right\}
			&\quad \text{(S)}\\
			\end{align*}

\paragraph{The Layer Problem at the right boundary} 			
			
			\begin{align*}
			&\left.\begin{array}{lll}
			\uuK = \sqrt{C} \quad &\Rightarrow\quad \uu(0) \in (-\sqrt{C},\infty),\ &\uuhat = C>0\\
			\uuK = 0 \quad &\Rightarrow\quad  \uu(0) \in [0,\infty),\ &\uuhat = 0
			\end{array}\right\}
			&\quad \text{(S) }\\
			&\quad \uuK = -\sqrt{C} \quad \Rightarrow\quad \uu_0=-\sqrt{C},\ \uuhat=C> 0
			&\quad \text{(U)}\\
			\end{align*}

Here $U$ denotes the unstable fixpoint and $S$ the stable one.

\subsection{Nodes with 2 edges}
\label{1-1}
We consider the case where arc i  is oriented into the node, arc j out of the node. Then, the kinetic coupling conditions  for \eqref{bgk} are
\begin{align}\label{eq:CC_11f}
f_{1}^i &=  f_{1}^j\\ f_{2}^j &=  f_{2}^i\ .
\end{align}
In this case,  one expects the coupled macroscopic solution to be given by the  solution of the Riemann problem for the Burgers equation at the coupling position without a kinetic layer.
Nevertheless, we go through the full procedure to explain it in the present simple case. 
As a first step we consider the kinetic $i$- and $j$-layers. 
\subsubsection{Combination of  layer analysis and  coupling conditions}
Expressing \eqref{eq:CC_11f} in the variables $\uu$ and $\uuhat$ we obtain
\begin{align*}
	\uu^i &= \uu^j\\
	\uuhat^i &= \uuhat^j\ .
\end{align*}
The two kinetic layers with asymptotic states $u_K^i$ and $u_K^j$ can be combined via the coupling conditions using the combinations discussed below.

\paragraph{{\bf Case 1, U-U}}
 $\uu^i = -\sqrt{C^i}$, $\uuhat^i=C^i$ with $C^i>0$ and  $\uu^j = \sqrt{C^j}$, $\uuhat^j=C^j$ with $C^j>0$.
	No solution is possible since $-\sqrt{C^i}\neq \sqrt{C^j}$.

\paragraph{{\bf Case 2, U-S}}	
	 $\uu^i = -\sqrt{C^i}$, $\uuhat^i=C^i$ with $C^i>0$ and  $\uu^j = \uu_0^j$, $\uuhat^j=C^j$ where $C^j\geq 0$ is fixed.
	Thus $\uuhat^i =\uuhat^j=C^j$ and $\uu_0^j=-\sqrt{C^i}$.
	\begin{align*}
		\uuK^j  = -\sqrt{C^j}\ .
	\end{align*}
\paragraph{{\bf Case 3, S-U}}	
	 $\uu^i = \uu_0^i$, $\uuhat^i=C^i$ where $C^i\geq 0$ is fixed and  $\uu^j = \sqrt{C^j}$, $\uuhat^j=C^j$ with $C^j>0$.
	Thus $\uuhat^j=\uuhat^i = C^i$ and $\uu^i_0=\sqrt{C^j}$.
	\begin{align*}
		\uuK^i = \sqrt{C^i}\ .
	\end{align*}
\paragraph{{\bf Case 4, S-S}}		
 $\uu^i = \uu_0^i$, $\uuhat^i=C^i$ where $C^i\geq 0$ is fixed and $\uu^j = \uu_0^j$, $\uuhat^j=C^j$ where $C^j\geq 0$ is fixed.
	Only possible if $C^i = C^j$, then $\uu_0^i = \uu_0^j\in (-\sqrt{C^i},\sqrt{C^i})$ free.
	\begin{align*}
	 \uuK^i =  \sqrt{C^i} \text{ and } \uuK^j = -\sqrt{C^j}, \;\text{ if } \;C^i = C^j
	\end{align*}

In a second step we combine the layer solutions with the solution of the Riemann problem between the asymptotic state of the layer $\uu_K$
and the value of the Burgers solution $\uu_B$.
The possible cases are discussed in the following.

\subsubsection{Combine Riemann Problem  with layer and coupling conditions}
\paragraph{{\bf Case 1, RP1-1}}
In this case  $\uu_B^i <0$ and $\uu_B^j >0$.
	\begin{align*}
		\uuK^i &\in(-\infty,0]:&(U) &\text{ or } ((S) \text{ with } C^i=0)\\
		 \uuK^j &\in [0,\infty):
	 & (U) &\text{ or } ((S) \text{ with } C^j=0)
	\end{align*}
	leads to the following subcases
	\begin{enumerate}
		\item[{\bf U-S}] with $C^j = 0$ . Thus $\uuK^i =\uuK^j =0$.
		\item [{\bf S-U}] with $C^i = 0$ . Thus $\uuK^i =\uuK^j =0$.
		\item[{\bf S-S}] with $C^i = C^j=0$ . Thus $\uuK^i =\uuK^j =0$. 
	\end{enumerate}
	Alltogether one obtains
	\begin{align*}
		 \uuK^i =\uuK^j =0\ .
	\end{align*}
\paragraph{{\bf Case2, RP1-2}} In this case $\uu_B^i <0$ and $\uu_B^j \leq 0$.
	\begin{align*}
	\uuK^i &\in (-\infty,0]:&(U) &\text{ or } ((S) \text{ with } C^i=0)\\
	 \uuK^j &\in \{\uu_B^j\}\cup [-\uu_B^j                   ,\infty):&
	 ((U) &\text{ with } \sqrt{C^j}\geq -\uu_B^j) \text{ or } ((S) \text{ with } \uu_B^j=-\sqrt{C^j})
	\end{align*}
	leads to the following subcases
	\begin{enumerate}
		\item[{\bf U-S}] with $C^j = (\uu_B^j)^2$ . Thus $\uuK^i =\uuK^j =\uu_B^j$.
		\item[{\bf S-U}] with $C^i = 0$ but $C^j\geq (\uu_B^j)^2$. Only if $\uu_B^j = 0$, then $\uuK^i =\uuK^i =0$.
		\item[{\bf S-S}] with $C^i = 0$, $\ C^j= (\uu_B^j)^2$. Only if $\uu_B^j = 0$, then $\uuK^i =\uuK^i =0$.
	\end{enumerate}
		Alltogether one obtains
	\begin{align*}
	 \uuK^i =\uuK^j =\uu_B^j\ .
	\end{align*}
\paragraph{{\bf Case 3, 	RP2-1}} Here $\uu_B^i > 0$ and $\uu_B^j > 0$.
	\begin{align*}
	\uuK^i &\in (-\infty,-\uu_B^i]\cup \{\uu_B^i\}: &
	((U) &\text{ with } \sqrt{C^i}\geq \uu_B^i) \text{ or } ((S) \text{ with } \uu_B^i=\sqrt{C^i}) \\
	 \uuK^j &\in  [0,\infty):
& (U) &\text{ or } ((S) \text{ with } C^j=0)
	\end{align*}
	leads to the following subcases
	\begin{enumerate}
		\item[{\bf U-S}] with $C^i\geq (\uu_B^i)^2$, but $C^j = 0$ . Only if $\uu_B^j = 0$, then $\uuK^i =\uuK^j =0$.
		\item[{\bf S-U}] with $C^i = (\uu_B^i)^2$. Thus $\uuK^i =\uuK^j =\uu_B^i$.
		\item[{\bf S-S}] with $\ C^i= (\uu_B^i)^2$, $C^j = 0$. Only if $\uu_B^j = 0$, then $\uuK^i =\uuK^i =0$.
	\end{enumerate}
		Alltogether one obtains
	\begin{align*}
	 \uuK^i =\uuK^j =\uu_B^i\ .
	\end{align*}
	\paragraph{{\bf Case 4, RP2-2}} Here $\uu_B^i > 0$ and $\uu_B^j \leq 0$.
	\begin{align*}
	\uuK^i &\in (-\infty,-\uu_B^i]\cup \{\uu_B^i\}:&
	((U) &\text{ with } \sqrt{C^i}\geq \uu_B^i) 
	\text{ or } ((S) \text{ with } \uu_B^i=\sqrt{C^i})\\
	  \uuK^j &\in \{\uu_B^j\}\cup [-\uu_B^j,\infty):
	 &  ((U) &\text{ with } \sqrt{C^j}\geq -\uu_B^j) 
	\text{ or } ((S) \text{ with } \uu_B^j=-\sqrt{C^j})
	\end{align*}
	leads to the following subcases
	\begin{enumerate}
		\item[{\bf U-S}] with $C^i\geq (\uu_B^i)^2$ and $C^j= (\uu_B^j)^2$ . Only if $\uu_B^i < -\uu_B^j$, 
		then $\uuK^j =\uuK^i =\uu_B^j$.
		\item[{\bf S-U}] with $C^i= (\uu_B^i)^2$ and $C^j\geq (\uu_B^j)^2$ . Only if $\uu_B^i >- \uu_B^j$, 
		then $\uuK^i =\uuK^j =\uu_B^i$.
		\item[{\bf S-S}] with $C^i= (\uu_B^i)^2$, $C^j= (\uu_B^j)^2$. Only if $\uu_B^i = -\uu_B^j$, 
		then $\uuK^i =\uuK^j =\uu_B^i=-\uu_B^j$.
	\end{enumerate}
		Alltogether one obtains
	\begin{align*}
	\begin{cases}
	\uuK^i =\uuK^j =\uu_B^i & \mbox{if} \; \uu_B^i\geq -\uu_B^j\\
	\uuK^i =\uuK^j =\uu_B^j & \mbox{if} \; \uu_B^i\leq -\uu_B^j\\
		\uuK^i =-\uuK^j =\uu_B^i=-\uu_B^j & \mbox{if} \; \uu_B^i = -\uu_B^j
	\end{cases} \ .
	\end{align*}
This means, as expected, no layers arise in this problem and  we  obtain exactly the solution of the Riemann problem for the Burgers equation at the junction.

\section{Coupling conditions for nodes with 3 edges}
\label{3node}

We have to distinguish between the different orientations of the nodes, i.e. we distinguish between
3-0,0-3,1-2 and 2-1 nodes, where the first number denotes the  number of ingoing edges and the second number the number of ougoing edges.

\subsection{3-0 and 0-3 nodes}
For the Burgers equation only outgoing or only ingoing edges are not interesting, since the macroscopic flux is always positive. 
Thus requesting conservation of mass
\begin{align*}
\sum_{i=1}^{3}(\rho^i)^2 = 0
\end{align*}
allows only for a trivial solution.

\subsection{1-2 node}
Arc $i$ is oriented into the node and arc $j$ and $k$ are oriented out of the node.
\begin{figure}[h]
	\begin{center}
	\externaltikz{sketch_12node}{
	\begin{tikzpicture}[thick]
	\def\len{2}
	\node[fill,circle] (N) at (0,0){};
	\draw[->] (-\len,0)--(N) node[above,pos = 0.5]{$i$};
	\draw[->] (0,0)--(\len,0.6) node[above,pos = 0.5]{$j$};
	\draw[->] (0,0)--(\len,-0.6) node[below,pos = 0.5]{$k$};
	\end{tikzpicture}
	}
	\end{center}
	\caption{1-2 node}
\end{figure}

We choose the symmetric kinetic coupling conditions similar to \cite{BKP16}
\begin{eqnarray}
f_1^i &= \frac{1}{2} (f_1^j+f_1^k)\\
f_2^j &= \frac{1}{2} (f_2^i+f_1^k)\\
f_2^k &= \frac{1}{2} (f_2^i+f_1^j)\ .
\end{eqnarray}
Note that these coupling conditions do only conserve the mass, if $v_1=-v_2$, which we assume from now on.
We reformulate the above conditions in terms of $\uu$ and $\uuhat$ assuming that $v_1$ and $v_2$ are identical on all edges
with $v_1=-v_2=-v$. One obtains
\begin{eqnarray}
\label{eq:CCBurgers12_rq_C1}
v\uu^i-\uuhat^i &= \frac{1}{2} \left(v\uu^j-\uuhat^j+v\uu^k-\uuhat^k\right)\\
\label{eq:CCBurgers12_rq_C2}
\uuhat^j+v \uu^j &= \frac{1}{2} \left(\uuhat^i+v\uu^i + v\uu^k-\uuhat^k\right)\\
\label{eq:CCBurgers12_rq_C3}
\uuhat^k+v \uu^k &= \frac{1}{2} \left(\uuhat^i+v\uu^i + v\uu^j-\uuhat^j\right)\ .
\end{eqnarray}

\begin{remark}
In general linear mass conserving kinetic coupling conditions in the 1-2 case have the form
\begin{eqnarray}
f_1^i =  \alpha f_1^j+ (1-\beta) f_1^k\\
f_2^j =  (1-\gamma) f_2^i+ \beta f_1^k\\
f_2^k = \gamma f_2^i+(1-\alpha)f_1^j
\end{eqnarray}
with 3 free parameters $\alpha,\beta,\gamma$.
\end{remark}

\subsubsection{Combine kinetic layer and coupling conditions}
First the combination of the coupling conditions  with the layer equations is considered.
We  face many different cases. 
In order to avoid repetitions we combine in this section the possible layer solutions with the coupling conditions.
Each layer can have either a stable solution (S) or an unstable solution (U).
Thus, for three edges we have eight combinations.

{\bf Case1, U-U-U.}	
	Inserting into the first equation of the coupling conditions \eqref{eq:CCBurgers12_rq_C1} we obtain
	\begin{align*}
		-v\sqrt{C^i}-C^i &= \frac{1}{2} \left(v\sqrt{C^j}-C^j+v\sqrt{C^k}-C^k\right)\\
		-v\sqrt{C^i}-\frac{C^i}{2} &= \frac{1}{2} \left(v\sqrt{C^j}+v\sqrt{C^k}\right)
		\end{align*}	
	with $C^i>0$, $C^j>0$ and $C^k>0$. The left hand side is negative, while the right hand side is positive.
	Thus this combination is not admissible.
	
{\bf Case 2, S-U-U} 
	Assume $C^i$ to be given.
	Inserting into the coupling conditions
	\begin{eqnarray*}
		v\uu_0^i-C^i &= \frac{1}{2} \left(v\sqrt{C^j}-C^j+v\sqrt{C^k}-C^k\right)\\
		C^j+v \sqrt{C^j} &= \frac{1}{2} \left(C^i+v\uu_0^i + v\sqrt{C^k}-C^k\right)\\
		C^k+v \sqrt{C^k} &= \frac{1}{2} \left(C^i+v\uu_0^i + v\sqrt{C^j}-C^j\right).
	\end{eqnarray*}	
	Combining one and three, as well as one and two, gives two identical expressions for $C^j$ and $C^k$. Thus $C^j=C^k$ as the problem is symmetric.
	Summing all three equations, we obtain the conservation of mass
	$	C^i=C^j+C^k$.
	Thus $C^j = C^k = \frac{1}{2}C^i$  and consequently $\uu_0^i=\frac{C^i}{2v}+\sqrt{\frac{C^i}{2}}$.
	
	{\bf Case3, U-S-U} (and {\bf U-U-S} analogously) 
	Assume $C^j\geq 0$ to be given. Then
	\begin{align*}
	-v\sqrt{C^i}-C^i &= \frac{1}{2} \left(v\uu_0^j-C^j+v\sqrt{C^k}-C^k\right)\\
	C^j+v \uu_0^j &= \frac{1}{2} \left(C^i-v\sqrt{C^i} + v\sqrt{C^k}-C^k\right)\\
	C^k+v \sqrt{C^k} &= \frac{1}{2} \left(C^i-v\sqrt{C^i} + v\uu_0^j-C^j\right).
	\end{align*}	
	Rearranging gives
	\begin{align*}
		C^k+\frac{1}{2}C^j+\frac{3v}{2}\left(\sqrt{C^k+C^j}+\sqrt{C^k}\right)=0.
	\end{align*}
	This equation does not allow for a solution $C^k>0$, since there are only positive terms on the left hand side.

	{\bf Case 4, S-S-U}(and {\bf S-U-S} analogously) Then $C^i,C^j\geq 0$ are given and  $C^k=C^i-C^j$. That means, there is only a solution, if $C^i> C^j$. The coupling conditions give
		\begin{align*}
			v\uu_0^i-C^i &= \frac{1}{2} \left(v\uu_0^j-C^j+v\sqrt{C^k}-C^k\right)\\
			C^j+v \uu_0^j &= \frac{1}{2} \left(C^i+v\uu_0^i + v\sqrt{C^k}-C^k\right)\\
			C^k+v \sqrt{C^k} &= \frac{1}{2} \left(C^i+v\uu_0^i + v\uu_0^j-C^j\right).
		\end{align*}	
		Solving  yields
		$\uu_0^i = \sqrt{C^k}+\frac{1}{3v}(C^k+C^i)$		
		and
		$\uu_0^j = \sqrt{C^k}+\frac{1}{3v}(C^k-C^j)$.
		Note that $\uu_0^j< \sqrt{C^j}$ should hold, which is only true for $C^k< C^j$.
		Combining this restriction with the conservation of mass we obtain $C^j< C^i< 2C^j$.
		
{\bf Case 5, U-S-S} $C^j,C^k\geq 0$ given. From the conservation of mass we obtain $C^i=C^j+C^k$.
		The corresponding values of $\uu_0^j$ and $\uu_0^j$ are obtained from
		\begin{align*}
			-v\sqrt{C^i}-C^i &= \frac{1}{2} \left(v\uu_0^j-C^j+v\uu_0^k-C^k\right)\\
			C^j+v \uu_0^j &= \frac{1}{2} \left(C^i-v\sqrt{C^i} + v\uu_0^k-C^k\right)\\
			C^k+v \uu_0^k &= \frac{1}{2} \left(C^i-v\sqrt{C^i} + v\uu_0^j-C^j\right).
		\end{align*}	
		We obtain
		\begin{align*}
			\uu_0^j = \frac{-1}{3v}(C^i+C^j)-\sqrt{C^i}\ ,
			&&
			\uu_0^k = \frac{-1}{3v}(C^i+C^k)-\sqrt{C^i}.
		\end{align*}
{\bf Case 6, S-S-S} $C^i,C^j,C^k\geq 0$ given.
		Then the conservation of mass requires $C^i=C^j+C^k$. Then 
		\begin{align*}
			v\uu_0^i-C^i &= \frac{1}{2} \left(v\uu_0^j-C^j+v\uu_0^k-C^k\right)\\
			C^j+v \uu_0^j &= \frac{1}{2} \left(C^i+v\uu_0^i + v\uu_0^k-C^k\right)\\
			C^k+v \uu_0^k &= \frac{1}{2} \left(C^i+v\uu_0^i + v\uu_0^j-C^j\right)\ ,
		\end{align*}	
		which is a linear system for $\uu_0^i,\uu_0^j,\uu_0^k$, but with rank $2$.
		The values of $\uu_0^i,\uu_0^j,\uu_0^k$ are not uniquely determined.
		
\subsubsection{Combine  kinetic layer  and Riemann problems}

Assuming  the  states $\uu_B^i,\uu_B^j,\uu_B^k$ to be given, we have to determine the new states $\uu_K^i,\uu_K^j$ and $\uu_K^k$ at the node.
We have to consider eight different configurations.
For  each of them all possible combinations with stable or unstable layer solutions have to be discussed.
Not admissible combinations are not listed.

{\bf Case 1, RP1-1-1}
	\begin{align*}
	\uuK^i &\in (-\infty,0] :
	&  (U) &\text{ or } ((S) \text{ with } C^i=0) \\
	\uuK^j &\in [0,\infty): 
	& (U) &\text{ or } ((S) \text{ with } C^j=0)\\
	 \uuK^k &\in [0,\infty):
	& (U) &\text{ or } ((S) \text{ with } C^k=0)
	\end{align*}	
	\begin{enumerate}
		\item[{\bf SUU}] with $C^i = 0$ $\Rightarrow C^j=C^k=0$ and $\uu_0^i=0$.
		\item[{\bf SSU}] with $C^i = C^j = 0$ $\Rightarrow C^k=0$ and $\uu_0^i=\uu_0^j=0$.
		\item[{\bf SUS}] with $C^i = C^k = 0$ $\Rightarrow C^j=0$ and $\uu_0^i=\uu_0^k=0$.
		\item[{\bf USS}] with $C^j = C^k = 0$ $\Rightarrow C^i=0$ and $\uu_0^j=\uu_0^k=0$.
		\item[{\bf SSS}] with $C^i = C^j = C^k = 0$ $\Rightarrow$ $\uu_0^i=\uu_0^j=\uu_0^k=0$.		
	\end{enumerate}
	This gives
	\begin{align*}
		\uuK^i &= 0 &\uuK^j &= 0 &\uuK^k &= 0 \\
		\uu_0^i &= 0 &\uu_0^j &= 0 &\uu_0^k &= 0 \ .
	\end{align*}
	
{\bf Case 2, RP1-1-2}
	\begin{align*}
	\uuK^i &\in (-\infty,0] :& (U) &\text{ or } ((S) \text{ with } C^i=0) \\
	 \uuK^j &\in [0,\infty): & (U)&\text{ or } ((S) \text{ with } C^j=0)\\
 \uuK^k &\in \{\uu_B^k\}\cup [-\uu_B^k,\infty):
	& (U) &\text{ with } C^k\geq (\uu_B^k)^2 \text{ or } ((S) \text{ with } C^k=(\uu_B^k)^2)
	\end{align*}	
	\begin{enumerate}
		\item[{\bf SUU}] with $C^i = 0$ $\Rightarrow C^j=C^k=0$, which contradicts $C^k\geq (\uu_B^k)^2 $, i.e. no solution.
		\item[{\bf SSU}] with $C^i = C^j = 0$ $\Rightarrow C^k=0$, which contradicts $C^k\geq (\uu_B^k)^2 $, i.e. no solution.
		\item[{\bf SUS}] with $C^i = 0$, $C^k = (\uu_B^k)^2 $. No solution, since we need $C^i\geq C^k$.
		\item[{\bf USS}] with $C^j = 0$, $C^k = (\uu_B^k)^2$ $\Rightarrow C^i=(\uu_B^k)^2$ and $\uu_0^j=\frac{-1}{3v}(\uu_B^k)^2+\uu_B^k$ and  $\uu_0^k=\frac{-2}{3v}(\uu_B^k)^2+\uu_B^k$, since $\uu_B^k<0$.
		\item[{\bf SSS}] with $C^i = C^j = 0$ and $C^k = (\uu_B^k)^2$ contradicts the conservation of mass, i.e. no solution.		
	\end{enumerate}
	Alltogether
	\begin{align*}
		\uuK^i &= \uu_B^k<0 &\uuK^j &= 0 &\uuK^k &= \uu_B^k<0 \\
		\uu_0^i &= \uu_B^k &\uu_0^j &= \frac{-1}{3v}(\uu_B^k)^2+\uu_B^k &\uu_0^k &= \frac{-2}{3v}(\uu_B^k)^2+\uu_B^k \ .
	\end{align*}
	
{\bf Case 3,RP1-2-1} This is symmetric to {\bf RP1-1-2}. We obtain
	\begin{align*}
	\uuK^i &= \uu_B^j<0 &\uuK^j &= \uu_B^j<0 &\uuK^k &= 0 \\
	\uu_0^i &= \uu_B^j &\uu_0^j &= \frac{-2}{3v}(\uu_B^j)^2+\uu_B^j &\uu_0^k &= \frac{-1}{3v}(\uu_B^j)^2+\uu_B^j \ .
	\end{align*}
	
{\bf Case 4,	RP2-1-1}
	\begin{align*}
	\uuK^i &\in (-\infty,-\uu_B^i]\cup \{\uu_B^i\}:& ((U) &\text{ with } C^i\geq (\uu_B^i)^2\text{ or } ((S) \text{ with } C^i=(\uu_B^i)^2) \\
	 \uuK^j &\in [0,\infty): & (U)&\text{ or } ((S) \text{ with } C^j=0)\\
 \uuK^k &\in [0,\infty):
	& (U)&\text{ or } ((S) \text{ with } C^k=0)
	\end{align*}	
	\begin{enumerate}
		\item[{\bf SUU}] with $C^i = (\uu_B^i)^2$ $\Rightarrow C^j=C^k=\frac{1}{2}(\uu_B^i)^2$ and $\uu_0^i=\frac{1}{2v}(\uu_B^i)^2+\frac{1}{\sqrt{2}}\uu_B^i$
		\item[{\bf SSU}] with $C^i = (\uu_B^i)^2$, $C^j = 0$ $\Rightarrow C^k=(\uu_B^i)^2$, but $C^i< 2C^j=0$ is violated, i.e. no solution.
		\item[{\bf SUS}] with $C^i = (\uu_B^i)^2$, $C^k = 0 $ $\Rightarrow$ $C^j = (\uu_B^i)^2$, but $C^i< 2C^k=0$ is violated, i.e. no solution.
		\item[{\bf USS}] with $C^j = 0$, $C^k = 0$ $\Rightarrow C^i=0$ which violates $ C^i\geq (\uu_B^i)^2$, i.e. no solution.
		\item[{\bf SSS}] with $C^i = (\uu_B^i)^2$ and $C^j = C^k = 0$ contradicts the conservation of mass, i.e. no solution.		
	\end{enumerate}
	Alltogether, we obtain
	\begin{align*}
		\uuK^i &= \uu_B^i>0 &\uuK^j &= \frac{1}{\sqrt{2}}\uu_B^i>0 &\uuK^k &= \frac{1}{\sqrt{2}}\uu_B^i>0 \\
		\uu_0^i &= \frac{1}{2v}(\uu_B^i)^2+\frac{1}{\sqrt{2}}\uu_B^i &\uu_0^j &= \frac{1}{\sqrt{2}}\uu_B^i &\uu_0^k&= \frac{1}{\sqrt{2}}\uu_B^i \ .
	\end{align*}
	
{\bf Case 5, RP1-2-2}
	\begin{align*}
	\uuK^i &\in (-\infty,0] : &(U) &\text{ or } ((S) \text{ with } C^i=0) \\
	 \uuK^j &\in \{\uu_B^j\}\cup [-\uu_B^j,\infty):& (U)& \text{ with } C^j\geq (\uu_B^j)^2 \text{ or } ((S) \text{ with } C^j=(\uu_B^j)^2)\\
 \uuK^k &\in \{\uu_B^k\}\cup [-\uu_B^k,\infty):
	& (U)& \text{ with } C^k\geq (\uu_B^k)^2 \text{ or } ((S) \text{ with } C^k=(\uu_B^k)^2)
	\end{align*}	
	\begin{enumerate}
		\item[{\bf SUU}] with $C^i = 0$ $\Rightarrow C^j=C^k=0$, which contradicts $C^j\geq (\uu_B^j)^2 $ and  $C^k\geq (\uu_B^k)^2 $, i.e. no solution.
		\item[{\bf SSU}] with $C^i = 0$, $C^j = (\uu_B^j)^2$, which contradicts the conservation of mass, i.e. no solution.
		\item[{\bf SUS}] with $C^i = 0$, $C^k = (\uu_B^k)^2 $, which contradicts the conservation of mass, i.e. no solution.
		\item[{\bf USS}] with $C^j = (\uu_B^j)^2$, $C^k = (\uu_B^k)^2$ $\Rightarrow C^i=(\uu_B^j)^2+(\uu_B^k)^2$ 
		and
		\begin{align*}
		\uu_0^j&=\frac{-1}{3v}\left(2(\uu_B^j)^2+(\uu_B^k)^2\right)-\sqrt{(\uu_B^2)^2+(\uu_B^k)^2}\\
		\uu_0^k&=\frac{-1}{3v}\left((\uu_B^j)^2+2(\uu_B^k)^2\right)-\sqrt{(\uu_B^2)^2+(\uu_B^k)^2}\ .
		\end{align*} 
				
		\item[{\bf SSS}] with $C^i = 0$, $C^j = (\uu_B^j)^2$ and $C^k = (\uu_B^k)^2$ contradicts the conservation of mass, i.e. no solution.		
	\end{enumerate}
	Alltogether, this is 
	\begin{align*}
	\uuK^i &= -\sqrt{(\uu_B^2)^2+(\uu_B^k)^2}<0 &\uu_0^i &= -\sqrt{(\uu_B^2)^2+(\uu_B^k)^2}\\
	\uuK^j &= \uu_B^j<0 &\uu_0^j &= \frac{-1}{3v}\left(2(\uu_B^j)^2+(\uu_B^k)^2\right)-\sqrt{(\uu_B^2)^2+(\uu_B^k)^2} \\
	\uuK^k &= \uu_B^k<0 & \uu_0^k &= \frac{-1}{3v}\left((\uu_B^j)^2+2(\uu_B^k)^2\right)-\sqrt{(\uu_B^2)^2+(\uu_B^k)^2} \ .
	\end{align*}

{\bf Case 6, RP2-1-2}
	\begin{align*}
	\uuK^i &\in (-\infty,-\uu_B^i]\cup \{\uu_B^i\} :&	((U) &\text{ with } C^i\geq (\uu_B^i)^2\text{ or } 	((S) \text{ with } C^i=(\uu_B^i)^2)\\
	 \uuK^j &\in [0,\infty) :& (U)&\text{ or } ((S) \text{ with } C^j=0)\\
  \uuK^k &\in \{\uu_B^k\}\cup [-\uu_B^k,\infty):&
	(U) &\text{ with } C^k\geq (\uu_B^k)^2 \text{ or } ((S) \text{ with } C^k=(\uu_B^k)^2)
	\end{align*}	
	\begin{enumerate}
		\item[{\bf SUU}] with $C^i = (\uu_B^i)^2$ $\Rightarrow C^j=C^k=\frac{1}{2}(\uu_B^i)^2$. Only possible, if $\frac{1}{2}(\uu_B^i)^2\geq (\uu_B^k)^2$.
		Then $\uu_0^i=\frac{1}{2v}(\uu_B^i)^2+\frac{1}{\sqrt{2}}\uu_B^i$.
		\item[{\bf SSU}] with $C^i = (\uu_B^i)^2$, $C^j = 0$ $\Rightarrow C^k=(\uu_B^i)^2$, but $C^i< 2C^j=0$ is violated, i.e. no solution.
		\item[{\bf SUS}] with $C^i = (\uu_B^i)^2$, $C^k = (\uu_B^k)^2 $ $\Rightarrow$ $C^j = (\uu_B^i)^2-(\uu_B^k)^2$, only if $(\uu_B^i)^2\geq (\uu_B^k)^2$ and $(\uu_B^i)^2< 2(\uu_B^k)^2$.
		$\uu_0^i = \frac{1}{3v}\left(2(\uu_B^i)^2-(\uu_B^k)^2\right)+\sqrt{(\uu_B^i)^2-(\uu_B^k)^2}$
		and $\uu_0^k = \frac{1}{3v}\left((\uu_B^i)^2-2(\uu_B^k)^2\right)+\sqrt{(\uu_B^i)^2-(\uu_B^k)^2}$.
		\item[{\bf USS}] with $C^j = 0$, $C^k = (\uu_B^k)^2$ $\Rightarrow C^i=(\uu_B^k)^2$ only if $ (\uu_B^k)^2> (\uu_B^i)^2$.
		$\uu_0^j = \frac{-1}{3v}\left((\uu_B^k)^2\right)+\uu_B^k$
		and $\uu_0^k = \frac{-2}{3v}\left((\uu_B^k)^2\right)+\uu_B^k$, since $\uu_B^k<0$.
		\item[{\bf SSS}] with $C^i = (\uu_B^i)^2$, $C^j = 0$ and $C^k = (\uu_B^k)^2$ is only admissible if $(\uu_B^i)^2=(\uu_B^k)^2$.
		Then $\uu_0^i=\uu_0^j=\uu_0^k=0$.
	\end{enumerate}
	This leads to the following three cases.
	
	If $\uu_B^i\geq -\sqrt{2}\uu_B^k$, then
		\begin{align*}
		&\uuK^i = \uu_B^i>0 &	\uu_0^i = \frac{1}{2v}(\uu_B^i)^2+\frac{1}{\sqrt{2}}\uu_B^i \\
		&\uuK^j = \frac{1}{\sqrt{2}}\uu_B^i>0& \uu_0^j = \frac{1}{\sqrt{2}}\uu_B^i\\
		&\uuK^k = \frac{1}{\sqrt{2}}\uu_B^i>0 &\uu_0^k= \frac{1}{\sqrt{2}}\uu_B^i \ ,
		\end{align*}
	if $-\sqrt{2}\uu_B^k> \uu_B^i\geq -\uu_B^k$, then
		\begin{align*}
		\uuK^i &= \uu_B^i>0 &\uu_0^i &= \frac{1}{3v}\left(2(\uu_B^i)^2-(\uu_B^k)^2\right)+\sqrt{(\uu_B^i)^2-(\uu_B^k)^2} \\
		\uuK^j &= \sqrt{(\uu_B^i)^2-(\uu_B^k)^2}>0 &\uu_0^j &= \sqrt{(\uu_B^i)^2-(\uu_B^k)^2} 
		\\
		\uuK^k &= \uu_B^k<0 &\uu_0^k&= \frac{1}{3v}\left((\uu_B^i)^2-2(\uu_B^k)^2\right)+\sqrt{(\uu_B^i)^2-(\uu_B^k)^2}\ ,
		\end{align*}
	if $\uu_B^i< -\uu_B^k$, then
		\begin{align*}
		\uuK^i &= \uu_B^k<0 &\uu_0^i &= \uu_B^k &\\
		\uuK^j &= 0 &\uu_0^j &= \frac{-1}{3v}\left((\uu_B^k)^2\right)+\uu_B^k\\
		\uuK^k &= \uu_B^k<0 &\uu_0^k&= \frac{-2}{3v}\left((\uu_B^k)^2\right)+\uu_B^k\ .
		\end{align*}

{\bf Case 7,	RP2-2-1}
	\begin{align*}
		\uuK^i &\in (-\infty,-\uu_B^i]\cup \{\uu_B^i\} :&	((U) &\text{ with } C^i\geq (\uu_B^i)^2 \text{ or } ((S) \text{ with } C^i=(\uu_B^i)^2) \\
		  \uuK^j &\in \{\uu_B^j\}\cup [-\uu_B^j,\infty):&
	 (U) &\text{ with } C^j\geq (\uu_B^j)^2 \text{ or } ((S) \text{ with } C^j=(\uu_B^j)^2)\\
		\uuK^k &\in [0,\infty): &
		(U)&\text{ or } ((S) \text{ with } C^k=0)
	\end{align*}	
	\begin{enumerate}
		\item[{\bf SUU}] with $C^i = (\uu_B^i)^2$ $\Rightarrow C^j=C^k=\frac{1}{2}(\uu_B^i)^2$. Only possible, if $\frac{1}{2}(\uu_B^i)^2\geq (\uu_B^j)^2$.
		Then $\uu_0^i=\frac{1}{2v}(\uu_B^i)^2+\frac{1}{\sqrt{2}}\uu_B^i$.
		\item[{\bf SSU}] with $C^i = (\uu_B^i)^2$, $C^j = (\uu_B^j)^2 $ $\Rightarrow$ $C^k = (\uu_B^i)^2-(\uu_B^j)^2$, only if $(\uu_B^i)^2\geq (\uu_B^j)^2$ and $(\uu_B^i)^2< 2(\uu_B^j)^2$.
		$\uu_0^i = \frac{1}{3v}\left(2(\uu_B^i)^2-(\uu_B^j)^2\right)+\sqrt{(\uu_B^i)^2-(\uu_B^j)^2}$
		and $\uu_0^k = \frac{1}{3v}\left((\uu_B^i)^2-2(\uu_B^k)^2\right)+\sqrt{(\uu_B^i)^2-(\uu_B^k)^2}$.
		\item[{\bf SUS}] with $C^i = (\uu_B^i)^2$, $C^k = 0$ $\Rightarrow C^j=(\uu_B^i)^2$, but $C^i< 2C^k=0$ is violated, i.e. no solution.
		\item[{\bf USS}] with $C^j = (\uu_B^j)^2$, $C^k = 0$ $\Rightarrow C^i=(\uu_B^j)^2$ only if $ (\uu_B^j)^2> (\uu_B^i)^2$.
		$\uu_0^j = \frac{-2}{3v}\left((\uu_B^j)^2\right)+\uu_B^j$
		and $\uu_0^k = \frac{-1}{3v}\left((\uu_B^j)^2\right)+\uu_B^j$, since $\uu_B^j<0$.
		\item[{\bf SSS}] with $C^i = (\uu_B^i)^2$, $C^j = (\uu_B^j)^2$ and $C^k = 0$ is only admissible if $(\uu_B^i)^2=(\uu_B^j)^2$.
		Then $\uu_0^i=\uu_0^j=\uu_0^k=0$.
	\end{enumerate}
	This leads to the following three cases.
	
	If $\uu_B^i\geq -\sqrt{2}\uu_B^j$, then
		\begin{align*}
		\uuK^i &= \uu_B^i>0 &\uu_0^i &= \frac{1}{2v}(\uu_B^i)^2+\frac{1}{\sqrt{2}}\uu_B^i\\
		\uuK^j &= \frac{1}{\sqrt{2}}\uu_B^i>0 &\uu_0^j &= \frac{1}{\sqrt{2}}\uu_B^i\\
		\uuK^k &= \frac{1}{\sqrt{2}}\uu_B^i>0 & \uu_0^k&= \frac{1}{\sqrt{2}}\uu_B^i \ ,
			\end{align*}
	if $-\sqrt{2}\uu_B^j> \uu_B^i\geq -\uu_B^j$
		\begin{align*}	
		\uuK^i &= \uu_B^i>0 &\uu_0^i &= \frac{1}{3v}\left(2(\uu_B^i)^2-(\uu_B^j)^2\right)+\sqrt{(\uu_B^i)^2-(\uu_B^j)^2} \\
		\uuK^j &= \uu_B^j<0  &\uu_0^j&= \frac{1}{3v}\left((\uu_B^i)^2-2(\uu_B^j)^2\right)+\sqrt{(\uu_B^i)^2-(\uu_B^j)^2}\\
		\uuK^k &= \sqrt{(\uu_B^i)^2-(\uu_B^j)^2}>0 &\uu_0^k &= \sqrt{(\uu_B^i)^2-(\uu_B^j)^2}
		\ ,
			\end{align*}
	if $\uu_B^i< -\uu_B^j$, then
		\begin{align*}	
		\uuK^i &= \uu_B^j<0 &\uu_0^i &= \uu_B^j\\
		\uuK^j &= \uu_B^j<0 & \uu_0^j&= \frac{-2}{3v}\left((\uu_B^j)^2\right)+\uu_B^j\\
		\uuK^k &= 0& \uu_0^k &= \frac{-1}{3v}\left((\uu_B^j)^2\right)+\uu_B^j \ .
		\end{align*}

{\bf Case 8, RP2-2-2}
	\begin{align*}
	\uuK^i &\in (-\infty,-\uu_B^i]\cup \{\uu_B^i\} :&	((U) &\text{ with } C^i\geq (\uu_B^i)^2\text{ or } ((S) \text{ with } C^i=(\uu_B^i)^2) \\
	 \uuK^j &\in \{\uu_B^j\}\cup [-\uu_B^j,\infty):
& (U)& \text{ with } C^j\geq (\uu_B^j)^2 \text{ or } ((S) \text{ with } C^j=(\uu_B^j)^2)\\
 \uuK^k &\in \{\uu_B^k\}\cup [-\uu_B^k,\infty):&
	(U) &\text{ with } C^k\geq (\uu_B^k)^2 \text{ or } ((S) \text{ with } C^k=(\uu_B^k)^2)
	\end{align*}	
	\begin{enumerate}
		\item[{\bf SUU}] with $C^i = (\uu_B^i)^2$ $\Rightarrow C^j=C^k=\frac{1}{2}(\uu_B^i)^2$. Only possible, if $\frac{1}{2}(\uu_B^i)^2> (\uu_B^j)^2$ and $\frac{1}{2}(\uu_B^i)^2> (\uu_B^k)^2$.
		Then $\uu_0^i=\frac{1}{2v}(\uu_B^i)^2+\frac{1}{\sqrt{2}}\uu_B^i$.
		\item[{\bf SSU}] with $C^i = (\uu_B^i)^2$, $C^j = (\uu_B^j)^2$ $\Rightarrow C^k=(\uu_B^i)^2-(\uu_B^j)^2$ only if $(\uu_B^j)^2+(\uu_B^k)^2< (\uu_B^i)^2< 2(\uu_B^j)^2$.
		Then $\uu_0^i = \uu_B^i+\frac{1}{3v}\left(2(\uu_B^i)^2-(\uu_B^j)^2\right)$
		and $\uu_0^i = \uu_B^i+\frac{1}{3v}\left((\uu_B^i)^2-2(\uu_B^j)^2\right)$.
		\item[{\bf SUS}] with $C^i = (\uu_B^i)^2$, $C^k = (\uu_B^k)^2 $ $\Rightarrow$ $C^j = (\uu_B^i)^2-(\uu_B^k)^2$, only if $(\uu_B^i)^2> (\uu_B^k)^2+(\uu_B^j)^2$ and $(\uu_B^i)^2< 2(\uu_B^k)^2$.
		Then $\uu_0^i = \frac{1}{3v}\left(2(\uu_B^i)^2-(\uu_B^k)^2\right)+\sqrt{(\uu_B^i)^2-(\uu_B^k)^2}$
		and $\uu_0^k = \frac{1}{3v}\left((\uu_B^i)^2-2(\uu_B^k)^2\right)+\sqrt{(\uu_B^i)^2-(\uu_B^k)^2}$.
		\item[{\bf USS}] with $C^j = (\uu_B^j)^2$, $C^k = (\uu_B^k)^2$ $\Rightarrow C^i=(\uu_B^j)^2+(\uu_B^k)^2$ only if $ (\uu_B^j)^2+(\uu_B^k)^2\geq (\uu_B^i)^2$.
		Then $\uu_0^j = \frac{-1}{3v}\left(2(\uu_B^j)^2+(\uu_B^k)^2\right)-\sqrt{(\uu_B^j)^2+(\uu_B^k)^2}$
		and $\uu_0^k = \frac{-1}{3v}\left((\uu_B^j)^2+2(\uu_B^k)^2\right)-\sqrt{(\uu_B^j)^2+(\uu_B^k)^2}$.
		\item[{\bf SSS}] with $C^i = (\uu_B^i)^2$, $C^j = (\uu_B^j)^2$ and $C^k = (\uu_B^k)^2$ is only admissible if $(\uu_B^i)^2=(\uu_B^j)^2+(\uu_B^k)^2$.
		Then $\uu_0^i=\uu_0^j=\uu_0^k=0$.
\end{enumerate}

	We have the following cases.

	If $\uu_B^i\geq -\sqrt{2}\uu_B^j, \uu_B^i\geq -\sqrt{2}\uu_B^k$, then	
		\begin{align*}
		\uuK^i &= \uu_B^i>0 &\uu_0^i &= \frac{1}{2v}(\uu_B^i)^2+\frac{1}{\sqrt{2}}\uu_B^i\\
		\uuK^j &= \frac{1}{\sqrt{2}}\uu_B^i>0 &\uu_0^j &= \frac{1}{\sqrt{2}}\uu_B^i\\
		\uuK^k &= \frac{1}{\sqrt{2}}\uu_B^i>0	&\uu_0^k&= \frac{1}{\sqrt{2}}\uu_B^i \ ,
		\end{align*}
	if $-\sqrt{2}\uu_B^j> \uu_B^i\geq \sqrt{(\uu_B^j)^2+(\uu_B^k)^2}$
		\begin{align*}	
		\uuK^i &= \uu_B^i>0 &\uu_0^i &= \frac{1}{3v}\left(2(\uu_B^i)^2-(\uu_B^j)^2\right)+\sqrt{(\uu_B^i)^2-(\uu_B^j)^2}\\
		\uuK^j &= \uu_B^j<0 &\uu_0^j&= \frac{1}{3v}\left((\uu_B^i)^2-2(\uu_B^j)^2\right)+\sqrt{(\uu_B^i)^2-(\uu_B^j)^2} \\
		\uuK^k &= \sqrt{(\uu_B^i)^2-(\uu_B^j)^2}>0 &\uu_0^k &= \sqrt{(\uu_B^i)^2-(\uu_B^j)^2}\ ,
		\end{align*}	

	if $-\sqrt{2}\uu_B^k> \uu_B^i\geq \sqrt{(\uu_B^j)^2+(\uu_B^k)^2}$
		\begin{align*}	
		\uuK^i &= \uu_B^i >0 &\uu_0^i &= \frac{1}{3v}\left(2(\uu_B^i)^2-(\uu_B^k)^2\right)+\sqrt{(\uu_B^i)^2-(\uu_B^k)^2}\\
		\uuK^j &= \sqrt{(\uu_B^i)^2-(\uu_B^k)^2}>0 & \uu_0^j &= \sqrt{(\uu_B^i)^2-(\uu_B^k)^2}>0\\
		\uuK^k &= \uu_B^k <0 &\uu_0^k &= \frac{1}{3v}\left((\uu_B^i)^2-2(\uu_B^k)^2\right)+\sqrt{(\uu_B^i)^2-(\uu_B^k)^2}\ ,
		\end{align*}	
	if $\uu_B^i\leq \sqrt{(\uu_B^j)^2+(\uu_B^k)^2}$
		\begin{align*}	
			\uuK^i &= -\sqrt{(\uu_B^j)^2+(\uu_B^k)^2}<0 &\uu_0^i &=  -\sqrt{(\uu_B^j)^2+(\uu_B^k)^2}<0\\
			\uuK^j &= \uu_B^j<0 &\uu_0^j &= \frac{-1}{3v}\left(2(\uu_B^j)^2+(\uu_B^k)^2\right)-\sqrt{(\uu_B^j)^2+(\uu_B^k)^2} \\
			\uuK^k &= \uu_B^k>0 &\uu_0^k &= \frac{-1}{3v}\left((\uu_B^j)^2+2(\uu_B^k)^2\right)-\sqrt{(\uu_B^j)^2+(\uu_B^k)^2}\ .
		\end{align*}	
		
	Note that these sub-cases partition uniquely the range of admissible states since ($x+y\geq z$) or ($x+y<z$ and $z<2x$) or ($x+y<z$ and $z<2y$) or ($z>2x$ and $z>2y$).

\subsubsection{Summary for nodes 1-2}

Alltogether one obtains the macroscopic coupling conditions

{\bf Case 1, RP1-1-1.} $\uu_B^i < 0 , \uu_B^j > 0 , \uu_B^k > 0$.
Then $C_i=C_j=C_k=0$ and
	\begin{align*}
		\uuK^i &= 0 &\uuK^j &= 0 &\uuK^k = 0 
	\end{align*}
	
{\bf Case 2, RP1-1-2.} $\uu_B^i <0, \uu_B^j>0, \uu_B^k<0$. Then  $
	C_i=C_k , C_j=0$ and	
	\begin{align*}
		\uuK^i &= \uu_B^k <0 &\uuK^j &= 0 &\uuK^k = \uu_B^k<0
	\end{align*}
	
{\bf Case 3, RP1-2-1} $\uu_B^i <0, \uu_B^j<0, \uu_B^k>0 $. Then
	$C_i=C_j,  C_k=0$ and 
	\begin{align*}
	\uuK^i &= \uu_B^j <0 &\uuK^j &= \uu_B^j<0 &\uuK^k = 0 
	\end{align*}

{\bf Case 4,	RP2-1-1} $\uu_B^i >0, \uu_B^j>0,\uu_B^k>0$. Then
	$C_j=\frac{C_i}{2},  C_k=\frac{C_i}{2}$ and 	
	\begin{align*}
		\uuK^i &= \uu_B^i>0 &\uuK^j &= \frac{1}{\sqrt{2}}\uu_B^i>0 &\uuK^k = \frac{1}{\sqrt{2}}\uu_B^i>0 
	\end{align*}
	
	{\bf Case 5, RP1-2-2} $\uu_B^i <0, \uu_B^j<0,\uu_B^k<0 $. Then 
		$C_i= C_j+C_k$ and 
	\begin{align*}
	\uuK^i &= -\sqrt{(\uu_B^2)^2+(\uu_B^k)^2}<0 &
	\uuK^j &= \uu_B^j<0 &
	\uuK^k = \uu_B^k<0 
	\end{align*}

{\bf Case 6, RP2-1-2}
		$\uu_B^i >0, \uu_B^j>0,\uu_B^k<0  $. Then we have 3 subcases.
	
	If $\uu_B^i\geq -\sqrt{2}\uu_B^k$, then $C_j=\frac{C_i}{2},  C_k=\frac{C_i}{2}$ and 
		\begin{align*}
		\uuK^i &= \uu_B^i>0 & 
		\uuK^j &= \frac{1}{\sqrt{2}}\uu_B^i>0&
		\uuK^k &= \frac{1}{\sqrt{2}}\uu_B^i>0 
		\end{align*}
		
		If 
		$-\sqrt{2}\uu_B^k> \uu_B^i\geq -\uu_B^k$, then 
		$C_j=C_i- C_k$ and 
		\begin{align*}
		\uuK^i &= \uu_B^i>0 &
		\uuK^j &= \sqrt{(\uu_B^i)^2-(\uu_B^k)^2}>0 &
		\uuK^k &= \uu_B^k<0 
		\end{align*}
		
	If 	$\uu_B^i< -\uu_B^k$, then
		$C_i=C_k,   C_j =0$ and 
		\begin{align*}
		\uuK^i &= \uu_B^k<0 &
		\uuK^j &= 0 &
		\uuK^k &= \uu_B^k<0 
		\end{align*}

{\bf Case 7,	RP2-2-1}
	$\uu_B^i >0, \uu_B^j<0, \uu_B^k>0  $. We have again 3 subcases.
	
	If 	
		$\uu_B^i\geq -\sqrt{2}\uu_B^j$, then
			$C_j=\frac{C_i}{2},   C_k=\frac{C_i}{2}$ and
		\begin{align*}
		\uuK^i &= \uu_B^i>0 &
		\uuK^j &= \frac{1}{\sqrt{2}}\uu_B^i>0 &
		\uuK^k &= \frac{1}{\sqrt{2}}\uu_B^i>0
			\end{align*}
			
	If 
	  $-\sqrt{2}\uu_B^j> \uu_B^i\geq -\uu_B^j $ then
		$C_k=C_i-C_j$ and
		\begin{align*}
		\uuK^i &= \uu_B^i>0 &
		\uuK^j &= \uu_B^j<0  &
		\uuK^k &= \sqrt{(\uu_B^i)^2-(\uu_B^j)^2}>0
			\end{align*}
			
		If 	
		$\uu_B^i< -\uu_B^j$, then
		$C_i=C_j,  C_k =0$ and 
			\begin{align*}
		\uuK^i &= \uu_B^j<0 &
		\uuK^j &= \uu_B^j<0 & 
		\uuK^k &= 0
	\end{align*}

{\bf Case 8, RP2-2-2}
$	\uu_B^i >0, \uu_B^j<0, \uu_B^k<0 $. We have 4 subcases. 

	If $ \uu_B^i\geq -\sqrt{2}\uu_B^j, \uu_B^i\geq -\sqrt{2}\uu_B^k$ then
		$C_j=\frac{C_i}{2}, C_k = \frac{C_i}{2}$ and 
			\begin{align*}
		\uuK^i &= \uu_B^i>0 &
		\uuK^j &= \frac{1}{\sqrt{2}}\uu_B^i>0 &
		\uuK^k &= \frac{1}{\sqrt{2}}\uu_B^i>0	
			\end{align*}
			
	If 
		$-\sqrt{2}\uu_B^j> \uu_B^i\geq \sqrt{(\uu_B^j)^2+(\uu_B^k)^2}$, then
		$C_k = C_i-C_j$ and
		\begin{align*}
		\uuK^i &= \uu_B^i>0 &
		\uuK^j &= \uu_B^j<0 &
		\uuK^k &= \sqrt{(\uu_B^i)^2-(\uu_B^j)^2}>0 
	\end{align*}	
	
	If 
	$-\sqrt{2}\uu_B^k> \uu_B^i\geq \sqrt{(\uu_B^j)^2+(\uu_B^k)^2}$, then
		$C_j = C_i-C_k$ and
			\begin{align*}	
		\uuK^i &= \uu_B^i >0 &
		\uuK^j &= \sqrt{(\uu_B^i)^2-(\uu_B^k)^2}>0 & 
		\uuK^k &= \uu_B^k <0 
	\end{align*}	
	
	If 
	$ \uu_B^i\leq \sqrt{(\uu_B^j)^2+(\uu_B^k)^2}$, then
			$C_i = C_j+C_k$ and 
				\begin{align*}	
			\uuK^i &= -\sqrt{(\uu_B^j)^2+(\uu_B^k)^2}<0 &
			\uuK^j &= \uu_B^j<0 & 			\uuK^k &= \uu_B^k>0 
		\end{align*}

\subsection{2-1 node}\label{sec:21-node}

In this section we consider the remaining case of a node with 2 incoming and one outgoing edge, see Figure \ref{2-1}.
The incoming edges are labelled by $i,j$, the outgoing edge is $k$.
\begin{figure}[h]
	\begin{center}
\externaltikz{sketch_21node}{
\begin{tikzpicture}[thick]
	\def\len{2}
	\node[fill,circle] (N) at (0,0){};
	\draw[->] (0,0)--(\len,0) node[above,pos = 0.5]{$k$};
	\draw[->] (-\len,0.6)--(N) node[above,pos = 0.5]{$i$};
	\draw[->] (-\len,0-0.6)--(N) node[below,pos = 0.5]{$j$};
	\end{tikzpicture}
}
	\end{center}
		\caption{2-1 node}
		\label{2-1}
\end{figure}

Analogously to the 1-2 node we choose the  kinetic coupling conditions
\begin{eqnarray}
f_1^i &= \frac{1}{2} (f_2^j+f_1^k)\\
f_1^j &= \frac{1}{2} (f_2^i+f_1^k)\\
f_2^k &= \frac{1}{2} (f_2^i+f_2^j)\ .
\end{eqnarray}

The procedure follows the case of a 1-2 node.
We state only the results.

{\bf Case 1, RP1-1-1	} $
	\uu_B^i < 0, \uu_B^j < 0, \uu_B^k > 0$, then	
	$C_i=C_j=C_k=0$ and 
		\begin{align*}
		\uuK^i &= 0 &\uuK^j &= 0 &\uuK^k = 0 
	\end{align*}
	
{\bf Case 2, RP1-1-2} $
	\uu_B^i <0, \uu_B^j<0, \uu_B^k<0 $, then
	$C_i=\frac{C_k}{2}, C_j=\frac{C_k}{2}$ and
		\begin{align*}
			\uuK^i &= \frac{1}{\sqrt{2}}\uu_B^k
			<0 &\uuK^j &= \frac{1}{\sqrt{2}}\uu_B^k<0 &\uuK^k =\uu_B^k<0 		
	\end{align*}
	
{\bf Case 3, RP1-2-1} 
$
	\uu_B^i <0, \uu_B^j>0, \uu_B^k>0  $, then
	$C_j=C_k ,  C_i=0$ and 
	\begin{align*}
	\uuK^i &= 0 &\uuK^j &= \uu_B^j>0 &\uuK^k = \uu_B^j  >0
	\end{align*}
	
{\bf Case 4,	RP2-1-1}
$
	\uu_B^i >0, \uu_B^j<0, \uu_B^k>0  $ and 
		$C_i=C_k,  C_j=0$ and 
		\begin{align*}
			\uuK^i &= \uu_B^i >0 &\uuK^j &= 0 &\uuK^k = \uu_B^i>0
	\end{align*}
	
{\bf Case 5, RP1-2-2} $
	\uu_B^i <0, \uu_B^j>0, \uu_B^k<0 
	$, then we have 3 subcases.

	If $\uu_B^k\leq -\sqrt{2}\uu_B^j$, then
		$C_i = C_j = C_k/2$ and
		\begin{align*}
		\uuK^i &= \frac{1}{\sqrt{2}}\uu_B^k<0 &
		\uuK^j &= \frac{1}{\sqrt{2}}\uu_B^k<0 &
		\uuK^k &= \uu_B^k<0
		\end{align*}
						
	If
		$-\sqrt{2}\uu_B^j< \uu_B^k< -\uu_B^j$, then
		$C_i =C_k-C_j$ and 
		\begin{align*}
		\uuK^i &= -\sqrt{(\uu_B^k)^2-(\uu_B^j)^2}<0 &
		\uuK^j &= \uu_B^j>0  &
		\uuK^k &= \uu_B^k>0 
			\end{align*}
			
	If
		$ \uu_B^k> -\uu_B^j$, then
		$C_i=0, C_j=C_k$ and 
		\begin{align*}
		\uuK^i &= 0 &
		\uuK^j &= \uu_B^j>0 & 
		\uuK^k &= \uu_B^j>0  
		\end{align*}

{\bf Case 6, RP2-1-2}
	$
		\uu_B^i >0,  \uu_B^j<0, \uu_B^k<0  
	$. We have again 3 subcases.
		
		If $\uu_B^k\leq -\sqrt{2}\uu_B^i$, then
			$C_i=C_j=C_k/2$
			\begin{align*}
			&\uuK^i = \frac{1}{\sqrt{2}}\uu_B^k<0 &
			&\uuK^j = \frac{1}{\sqrt{2}}\uu_B^k<0& 
			&\uuK^k = \uu_B^k<0
			\end{align*}
			
			If 
			$-\sqrt{2}\uu_B^i<\uu_B^k< -\uu_B^i$, then
			$C_j = C_k -C_i$ and
			\begin{align*}
			\uuK^i &= \uu_B^i>0 &
			\uuK^j &= -\sqrt{(\uu_B^k)^2-(\uu_B^i)^2}<0 &
			\uuK^k = \uu_B^k<0 
			\end{align*}	
		
	If 
		$\uu_B^k< -\uu_B^i$, then
		$C_i=C_k, C_j =0$ and 
		\begin{align*}
		\uuK^i &= \uu_B^i>0 &
		\uuK^j &= 0 &
		\uuK^k &= \uu_B^i>0 
		\end{align*}
				
{\bf Case 7,	RP2-2-1}
		$\uu_B^i >0, \uu_B^j>0, \uu_B^k>0 $, then
		$C_k = C_i + C_j $ and
	\begin{align*}
	\uuK^i &= \uu_B^i>0 &
	\uuK^j &= \uu_B^j>0 &
	\uuK^k &= \sqrt{(\uu_B^i)^2+(\uu_B^j)^2}>0  
	\end{align*}

{\bf Case 8, RP2-2-2}
$
		\uu_B^i >0, \uu_B^j>0, \uu_B^k<0 $ . We have 4 subcases.
		
		If $ \uu_B^k\leq -\sqrt{2}\uu_B^i, \uu_B^k\leq -\sqrt{2}\uu_B^j$, then
		$C_j=\frac{C_i}{2}, C_k = \frac{C_i}{2}$ and
		\begin{align*}
		\uuK^i &= \frac{1}{\sqrt{2}}\uu_B^k>0 &
		\uuK^j &= \frac{1}{\sqrt{2}}\uu_B^k>0 &
		\uuK^k &= \uu_B^k<0	
			\end{align*}
			
			If
	$-\sqrt{2}\uu_B^j < \uu_B^k\leq -\sqrt{(\uu_B^i)^2+(\uu_B^j)^2}$, then
		$C_j = C_k-C_i$ and
		\begin{align*}
		\uuK^i &= \uu_B^i>0 &
		\uuK^j &= -\sqrt{(\uu_B^k)^2-(\uu_B^i)^2}<0  &
		\uuK^k &= \uu_B^k<0
	\end{align*}	
	
If 
$-\sqrt{2}\uu_B^i < \uu_B^k\leq -\sqrt{(\uu_B^i)^2+(\uu_B^j)^2}$, then
		$C_i = C_k-C_j$ and 
			\begin{align*}
		\uuK^i &= -\sqrt{(\uu_B^k)^2-(\uu_B^j)^2}<0 &
		\uuK^j &=  \uu_B^j >0& 
		\uuK^k &= \uu_B^k <0 
	\end{align*}	
	
	If 
		$ \uu_B^k\geq -\sqrt{(\uu_B^i)^2+(\uu_B^j)^2}$, then
			$C_i = C_j+C_k$ and
			\begin{align*}
			\uuK^i &= \uu_B^i>0 &
			\uuK^j &= \uu_B^j>0 &
			\uuK^k &= \sqrt{(\uu_B^i)^2+(\uu_B^j)^2}>0  
		\end{align*}

Note that all states and switching conditions are independent of $v$.

\section{Numerical results}
\label{Numerical results}
In this section the derived coupling conditions are investigated numerically.
The numerical solutions of the macroscopic equation \eqref{eq:CL} with $F(u) = u^2$ are compared to those obtained for the kinetic model \eqref{bgk}.

As numerical scheme for the kinetic equations the Upwind method for the linear advective part is combined with an implicit Euler scheme for the source term. 
The solution of the macroscopic equation is approximated with a Godunov scheme.
For all computations $1000$ cells are used as spacial resolution and the time steps are chosen according to the respective CFL conditions. 
The simulations are computed up to time $T=0.5$ and the relaxation parameter is chosen as $\epsilon = 0.0005$.
The initial conditions are formulated in macroscopic states, the remaining values in the kinetic model are chosen according to the relaxed state, i.e. $\hat{u}=F(u)$.
The kinetic speeds are chosen in agreement with the subcharacteristic condition as $v_1 = -2$ and $v_2=2$.
Larger values of $v$ would lead to similar results, as the coupling conditions are independent of $v$.
\subsection{Boundary conditions}

\begin{figure}[ht!]	
	\centering										
	\externaltikz{boundaryLayerplot}{
		\begin{tikzpicture}[scale=1.0]
		\begin{axis}[
		y=3cm,
		x=8cm,
		ylabel = $\uu$,xlabel = $x$,
		legend style = {at={(0.5,1)},xshift=0.2cm,yshift=0.1cm,anchor=south},
		legend columns= 3,
		]
		\addplot[color = blue,thick] file{Data/rho_kinetic_layertest1.txt};
		\addlegendentry{kinetic}
		\addplot[color = red,thick] file{Data/rho_Burgers_layertest1.txt};
		\addlegendentry{Burgers}
		\addplot[color = green,thick] file{Data/rho_layerL_layertest1.txt};
		\addplot[color = green,thick] file{Data/rho_layerR_layertest1.txt};
		\addlegendentry{Layer}
		\end{axis}
		\end{tikzpicture}
	}
	\caption{Boundary layer at the left and right end for $\epsilon = 0.0005$, $\uu_0(x) = 0.5$ and $f_2(0,t)= \frac{-1}{4}$, $f_1(1,t) = \frac{-9}{64}$.}
	\label{fig:BoundaryLayer}
\end{figure}
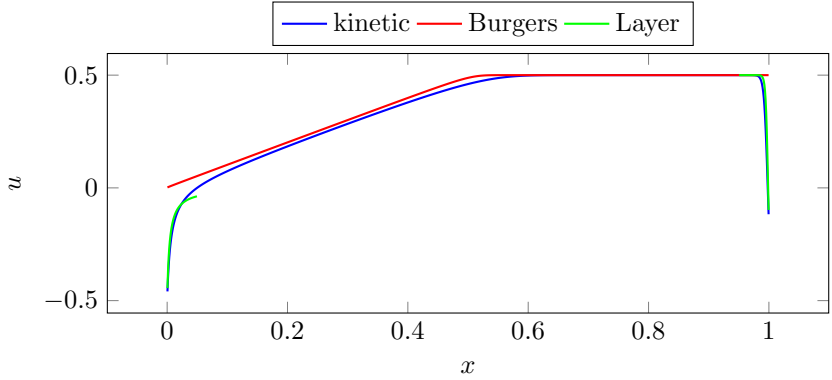

In this first example the boundary layer in the kinetic equation is  illustrated.
In Figure \ref{fig:BoundaryLayer} the numerical solution of the kinetic equation with the boundary values $f_2(0,t)= \frac{-1}{4}$ and $f_1(1,t) = \frac{-9}{64}$ is shown. 
The red line represents the solution of the Burgers equation with the corresponding boundary values obtained from section \ref{sec:MacroBoundaryCond}.
Additionally the solutions of the half-space problems are plotted (green lines), where the $x$-values are scaled with $\frac{1}{\epsilon}$.
On the left hand side the layer $\frac{-4}{x-\frac{4}{-0.5}}$ matches well the values at the boundary. 
With increasing distance to the boundary the layer differs from the kinetic solution since it can not take into account the rarefaction wave starting at $x=0$.
On the right boundary the layer is almost identical to the solution of the kinetic equation (blue line).

\subsection{Junctions}
For the junctions we consider the $1$-$2$-junction and the $2$-$1$-junction separately.
In both cases the coupling conditions are tested with six different Riemann problems at the junction. 

\subsubsection{Junction 1-2}

	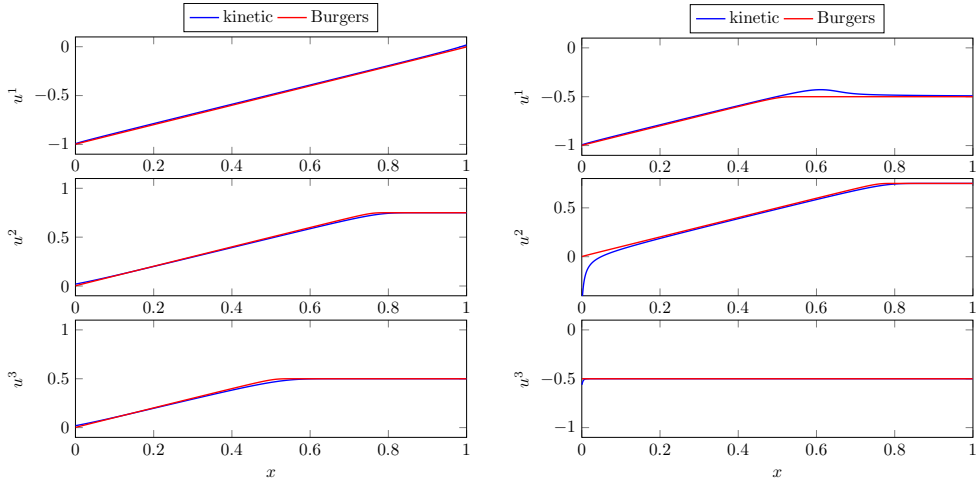
\begin{figure}
		\externaltikz{junction12_111-112}{
	\begin{tikzpicture}[scale=0.65]
		\begin{axis}[ylabel = $u^1$,
		y=2cm,
		x=8cm,
		xmin = -0.0, xmax = 1.0,
		ymin = -1.1, ymax = 0.1,
		legend style = {at={(0.5,1)},xshift=0.2cm,yshift=0.1cm,anchor=south},
		legend columns= 3,
		name=edge1,
		]
		\addplot[color = blue,thick] file{Data/rho_kinetic_junction12_test1_iE1.txt};
		\addlegendentry{kinetic}
		\addplot[color = red,thick] file{Data/rho_Burgers_junction12_test1_iE1.txt};
		\addlegendentry{Burgers}
		\end{axis}

		\begin{axis}[ylabel = $u^2$,
		y=2cm,
		x=8cm,
		xmin = -0.0, xmax = 1.0,
		ymin = -0.1, ymax = 1.1,
		 name=edge2,
		 at=(edge1.below south west), anchor=above north west,
		]
		\addplot[color = blue,thick] file{Data/rho_kinetic_junction12_test1_iE2.txt};
		\addplot[color = red,thick] file{Data/rho_Burgers_junction12_test1_iE2.txt};
		\end{axis}

		\begin{axis}[ylabel = $u^3$,xlabel =  $x$,
		y=2cm,
		x=8cm,
		xmin = -0.0, xmax = 1.0,
		ymin = -0.1, ymax = 1.1,
		 at=(edge2.below south west), anchor=above north west,
		]
		\addplot[color = blue,thick] file{Data/rho_kinetic_junction12_test1_iE3.txt};
		\addplot[color = red,thick] file{Data/rho_Burgers_junction12_test1_iE3.txt};
		\end{axis}
	\end{tikzpicture}
	\quad 
	\begin{tikzpicture}[scale=0.65]
		\begin{axis}[ylabel = $u^1$,
		y=2cm,
		x=8cm,
		xmin = -0.0, xmax = 1.0,
		ymin = -1.1, ymax = 0.1,
		legend style = {at={(0.5,1)},xshift=0.2cm,yshift=0.1cm,anchor=south},
		legend columns= 3,
		name=edge1,
		]
		\addplot[color = blue,thick] file{Data/rho_kinetic_junction12_test2_iE1.txt};
		\addlegendentry{kinetic}
		\addplot[color = red,thick] file{Data/rho_Burgers_junction12_test2_iE1.txt};
		\addlegendentry{Burgers}
		\end{axis}
		
		\begin{axis}[ylabel = $u^2$,
		y=2cm,
		x=8cm,
		xmin = -0.0, xmax = 1.0,
		ymin = -0.4, ymax = 0.8,
		name=edge2,
		at=(edge1.below south west), anchor=above north west,
		]
		\addplot[color = blue,thick] file{Data/rho_kinetic_junction12_test2_iE2.txt};
		\addplot[color = red,thick] file{Data/rho_Burgers_junction12_test2_iE2.txt};
		\end{axis}
		
		\begin{axis}[ylabel = $u^3$,xlabel =  $x$,
		y=2cm,
		x=8cm,
		xmin = -0.0, xmax = 1.0,
		ymin = -1.1, ymax = 0.1,
		at=(edge2.below south west), anchor=above north west,
		]
		\addplot[color = blue,thick] file{Data/rho_kinetic_junction12_test2_iE3.txt};
		\addplot[color = red,thick] file{Data/rho_Burgers_junction12_test2_iE3.txt};
		\end{axis}
	\end{tikzpicture}
	}
	\caption{Left: RP$1$-$1$-$1$ with $u^1=-1,$  $u^2=0.75,$  $u^3=0.5$, Right: RP$1$-$1$-$2$  with $u^1=-1,$  $u^2=0.75,$  $u^3=-0.5$}
		\label{fig:Network111_112}
	\end{figure}
	
In the $1$-$2$-junction tests six flow scenarios are investigated
covering  all  relevant states at the junction.
The first edge is connected to the junction at $x=1$, while the other two  edges are connected to the node at $x=0$.
Thus in  the following figures waves move to the left in the first edge but to the right in edge $2$ and $3$.
	
In Figure \ref{fig:Network111_112} the solutions to the initial conditions $(u_0^1,u_0^2,u_0^3)(x)=(-1, 0.75, 0.5)$ and $(u_0^1,u_0^2,u_0^3)(x)=(-1, 0.75, -0.5)$ are shown.
These correspond to the cases RP$1$-$1$-$1$ and RP$1$-$1$-$2$ respectively.
In RP$1$-$1$-$1$ the initial states only have characteristic speeds away from the junction and the coupling conditions enforce zero states in all edges. 
This leads to three rarefaction waves.  

On the right hand side RP$1$-$1$-$2$ is considered.
Flow is entering from edge $3$ but only exiting in edge $1$.
In edge $1$ a rarefaction wave forms and moves to the left.
On the slower end of the rarefaction wave a bump in the kinetic solution is present.
As the initial states at $t=0$ do not satisfy the coupling conditions, this small disturbance arises due to the transition in the new state at the junction.
For smaller values of $\epsilon$ and when refining the spacial and the temporal grid, this disturbance becomes narrower and more peaked.
Such temporal layers due to the initial conditions will also occur in other test cases.
On edge $2$ a boundary layer connects the junction state and a rarefaction wave, similar to the situation in Figure \ref{fig:BoundaryLayer}.
The ingoing flow from edge $3$ leads to a small boundary layer.

	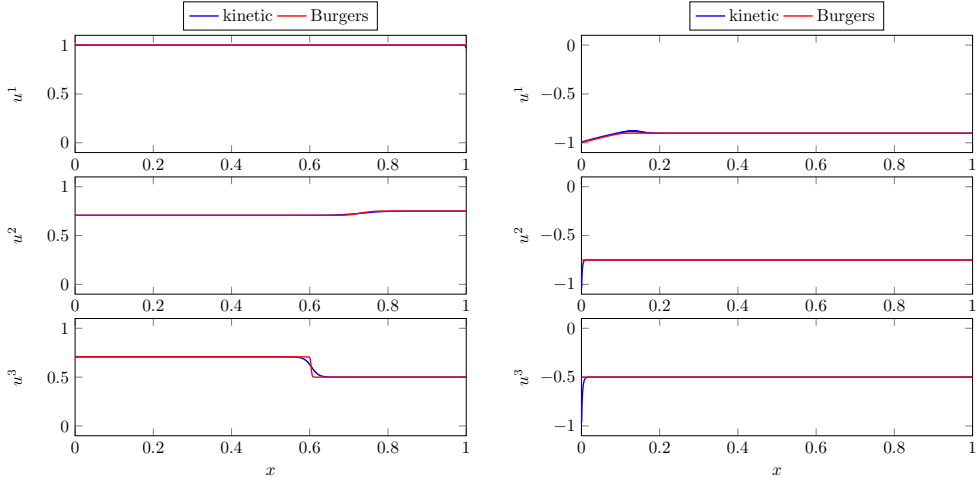
\begin{figure}		
		\externaltikz{junction12_211-122}{
		\begin{tikzpicture}[scale=0.65]
		\begin{axis}[ylabel = $u^1$,
		y=2cm,
		x=8cm,
		xmin = -0.0, xmax = 1.0,
		ymin = -0.1, ymax = 1.1,
		legend style = {at={(0.5,1)},xshift=0.2cm,yshift=0.1cm,anchor=south},
		legend columns= 3,
		name=edge1,
		]
		\addplot[color = blue,thick] file{Data/rho_kinetic_junction12_test3_iE1.txt};
		\addlegendentry{kinetic}
		\addplot[color = red,thick] file{Data/rho_Burgers_junction12_test3_iE1.txt};
		\addlegendentry{Burgers}
		\end{axis}
		
		\begin{axis}[ylabel = $u^2$,
		y=2cm,
		x=8cm,
		xmin = -0.0, xmax = 1.0,
		ymin = -0.1, ymax = 1.1,
		name=edge2,
		at=(edge1.below south west), anchor=above north west,
		]
		\addplot[color = blue,thick] file{Data/rho_kinetic_junction12_test3_iE2.txt};
		\addplot[color = red,thick] file{Data/rho_Burgers_junction12_test3_iE2.txt};
		\end{axis}
		
		\begin{axis}[ylabel = $u^3$,xlabel =  $x$,
		y=2cm,
		x=8cm,
		xmin = -0.0, xmax = 1.0,
		ymin = -0.1, ymax = 1.1,
		at=(edge2.below south west), anchor=above north west,
		]
		\addplot[color = blue,thick] file{Data/rho_kinetic_junction12_test3_iE3.txt};
		\addplot[color = red,thick] file{Data/rho_Burgers_junction12_test3_iE3.txt};
		\end{axis}
		\end{tikzpicture}
		\quad 
		\begin{tikzpicture}[scale=0.65]
		\begin{axis}[ylabel = $u^1$,
		y=2cm,
		x=8cm,
		xmin = -0.0, xmax = 1.0,
		ymin = -1.1, ymax = 0.1,
		legend style = {at={(0.5,1)},xshift=0.2cm,yshift=0.1cm,anchor=south},
		legend columns= 3,
		name=edge1,
		]
		\addplot[color = blue,thick] file{Data/rho_kinetic_junction12_test4_iE1.txt};
		\addlegendentry{kinetic}
		\addplot[color = red,thick] file{Data/rho_Burgers_junction12_test4_iE1.txt};
		\addlegendentry{Burgers}
		\end{axis}
		
		\begin{axis}[ylabel = $u^2$,
		y=2cm,
		x=8cm,
		xmin = -0.0, xmax = 1.0,
		ymin = -1.1, ymax = 0.1,
		name=edge2,
		at=(edge1.below south west), anchor=above north west,
		]
		\addplot[color = blue,thick] file{Data/rho_kinetic_junction12_test4_iE2.txt};
		\addplot[color = red,thick] file{Data/rho_Burgers_junction12_test4_iE2.txt};
		\end{axis}
		
		\begin{axis}[ylabel = $u^3$,xlabel =  $x$,
		y=2cm,
		x=8cm,
		xmin = -0.0, xmax = 1.0,
		ymin = -1.1, ymax = 0.1,
		at=(edge2.below south west), anchor=above north west,
		]
		\addplot[color = blue,thick] file{Data/rho_kinetic_junction12_test4_iE3.txt};
		\addplot[color = red,thick] file{Data/rho_Burgers_junction12_test4_iE3.txt};
		\end{axis}
		\end{tikzpicture}
		}
		\caption{Left: RP$2$-$1$-$1$  with $u^1=1,$  $u^2=0.75,$  $u^3=0.5$, Right: RP$1$-$2$-$2$ with $u^1=-1,$  $u^2=-0.75,$  $u^3=-0.5$}
		\label{fig:Network221_122}
	\end{figure}

	In Figure \ref{fig:Network221_122} on the left the flow from the first edge is split to the outgoing edges. 
	On the right hand side the flow from edge $2$ and $3$ is merged into edge $1$.
	Two layers connect the backward going flows with the junction states.
	In the first edge a small rarefaction wave travels to the left, followed by a small temporal layer. 
	This corresponds to the case$1$-$2$-$2$.
		
	\begin{figure}
		\externaltikz{junction12_212-222}{
			\begin{tikzpicture}[scale=0.65]
			\begin{axis}[ylabel = $u^1$,
			y=2cm,
			x=8cm,
		xmin = -0.0, xmax = 1.0,
			ymin = -0.1, ymax = 1.1,
			legend style = {at={(0.5,1)},xshift=0.2cm,yshift=0.1cm,anchor=south},
			legend columns= 3,
			name=edge1,
			]
			\addplot[color = blue,thick] file{Data/rho_kinetic_junction12_test5_iE1.txt};
			\addlegendentry{kinetic}
			\addplot[color = red,thick] file{Data/rho_Burgers_junction12_test5_iE1.txt};
			\addlegendentry{Burgers}
			\end{axis}
			
			\begin{axis}[ylabel = $u^2$,
			y=2cm,
			x=8cm,
			xmin = -0.0, xmax = 1.0,
			ymin = -.1, ymax = 1.1,
			name=edge2,
			at=(edge1.below south west), anchor=above north west,
			]
			\addplot[color = blue,thick] file{Data/rho_kinetic_junction12_test5_iE2.txt};
			\addplot[color = red,thick] file{Data/rho_Burgers_junction12_test5_iE2.txt};
			\end{axis}
			
			\begin{axis}[ylabel = $u^3$,xlabel =  $x$,
			y=2cm,
			x=8cm,
			xmin = -0.0, xmax = 1.0,
			ymin = -0.9, ymax = 0.3,
			at=(edge2.below south west), anchor=above north west,
			]
			\addplot[color = blue,thick] file{Data/rho_kinetic_junction12_test5_iE3.txt};
			\addplot[color = red,thick] file{Data/rho_Burgers_junction12_test5_iE3.txt};
			\end{axis}
			\end{tikzpicture}
			\quad 
			\begin{tikzpicture}[scale=0.65]
			\begin{axis}[ylabel = $u^1$,
			y=2cm,
			x=8cm,
			xmin = -0.0, xmax = 1.0,
			ymin = -1.1, ymax = 0.9,
			legend style = {at={(0.5,1)},xshift=0.2cm,yshift=0.1cm,anchor=south},
			legend columns= 3,
			name=edge1,
			]
			\addplot[color = blue,thick] file{Data/rho_kinetic_junction12_test6_iE1.txt};
			\addlegendentry{kinetic}
			\addplot[color = red,thick] file{Data/rho_Burgers_junction12_test6_iE1.txt};
			\addlegendentry{Burgers}
			\end{axis}
			
			\begin{axis}[ylabel = $u^2$,
			y=2cm,
			x=8cm,
			xmin = -0.0, xmax = 1.0,
			ymin = -1.1, ymax = -0.4,
			ytick distance=0.5,
			name=edge2,
			at=(edge1.below south west), anchor=above north west,
			]
			\addplot[color = blue,thick] file{Data/rho_kinetic_junction12_test6_iE2.txt};
			\addplot[color = red,thick] file{Data/rho_Burgers_junction12_test6_iE2.txt};
			\end{axis}
			
			\begin{axis}[ylabel = $u^3$,xlabel =  $x$,
			y=2cm,
			x=8cm,
			xmin = -0.0, xmax = 1.0,
			ymin = -1.1, ymax = -0.3,
			ytick distance=0.5,
			at=(edge2.below south west), anchor=above north west,
			]
			\addplot[color = blue,thick] file{Data/rho_kinetic_junction12_test6_iE3.txt};
			\addplot[color = red,thick] file{Data/rho_Burgers_junction12_test6_iE3.txt};
			\end{axis}
			\end{tikzpicture}
		}
		\caption{Left: RP$2$-$1$-$2$ with  $u^1=0.6,$  $u^2=0.75,$  $u^3=-0.5$, Right: RP$2$-$2$-$2$ with $u^1=0.8,$  $u^2=-0.75,$  $u^3=-0.5$}
		\label{fig:Network212_222}
	\end{figure}
	
	Figure \ref{fig:Network212_222} shows the last two test-cases for this  junction.
	On the left in the case RP$2$-$1$-$2$ the flow enters from the first and the third edge and exits into the second one, where a rarefaction wave moves to the right.
	The test on the right enforces states at the junctions similar to those of RP$1$-$2$-$2$ from Figure \ref{fig:Network221_122}.
	Here a strong transsonic rarefaction wave travels to the left in the first edge.
	In edges $2$ and $3$ two boundary layers form.
	Note that a linearization approach at the junction would fail in this test, as the characteristics of all initial states point into the node and no coupling conditions could be set. 
	The present coupling enforces the change of sign in edge $3$ in order to assure the conservation of mass.

	In all tests of the $1$-$2$ junction the kinetic and the macroscopic solutions are very close.
	Especially the states at the junction are correctly represented by the derived coupling conditions.
	Since the value of $\epsilon$ is small, also the boundary layers in the kinetic solution have a small spacial width.

\subsubsection{Junction 2-1}
	At the $2$-$1$ junction the edges $1$ and $2$ are orientated towards the coupling point, i.e. they are connected at $x=1$ to the junction.
	Edge $3$ is coupled at $x=0$.
	When studying Riemann Problems at the junction, the waves in the first two edges move to the left and in edge $3$ to the right.

\begin{figure}
	\externaltikz{junction21_111_112}{
		\begin{tikzpicture}[scale=0.65]
		\begin{axis}[ylabel = $u^1$,
		y=2cm,
		x=8cm,
		xmin = -0.0, xmax = 1.0,
		ymin = -1.1, ymax = .1,
		legend style = {at={(0.5,1)},xshift=0.2cm,yshift=0.1cm,anchor=south},
		legend columns= 3,
		name=edge1,
		]
		\addplot[color = blue,thick] file{Data/rho_kinetic_junction21_test1_iE1.txt};
		\addlegendentry{kinetic}
		\addplot[color = red,thick] file{Data/rho_Burgers_junction21_test1_iE1.txt};
		\addlegendentry{Burgers}
		\end{axis}
		
		\begin{axis}[ylabel = $u^2$,
		y=2cm,
		x=8cm,
		xmin = -0.0, xmax = 1.0,
		ymin = -1.1, ymax = .1,
		name=edge2,
		at=(edge1.below south west), anchor=above north west,
		]
		\addplot[color = blue,thick] file{Data/rho_kinetic_junction21_test1_iE2.txt};
		\addplot[color = red,thick] file{Data/rho_Burgers_junction21_test1_iE2.txt};
		\end{axis}
		
		\begin{axis}[ylabel = $u^3$,xlabel =  $x$,
		y=2cm,
		x=8cm,
		xmin = -0.0, xmax = 1.0,
		ymin = -0.1, ymax = 1.1,
		at=(edge2.below south west), anchor=above north west,
		]
		\addplot[color = blue,thick] file{Data/rho_kinetic_junction21_test1_iE3.txt};
		\addplot[color = red,thick] file{Data/rho_Burgers_junction21_test1_iE3.txt};
		\end{axis}
		\end{tikzpicture}
		\begin{tikzpicture}[scale=0.65]
		\begin{axis}[ylabel = $u^1$,
		y=2cm,
		y=2cm,
		width=8cm,
		xmin = -0.0, xmax = 1.0,
		scaled x ticks=false,
		ymin = -1.1, ymax = 0.1,
		legend style = {at={(0.5,1)},xshift=0.2cm,yshift=0.1cm,anchor=south},
		legend columns= 3,
		name=edge1,
		]
		\addplot[color = blue,thick] file{Data/rho_kinetic_junction21_test2_iE1.txt};
		\addlegendentry{kinetic}
		\addplot[color = red,thick] file{Data/rho_Burgers_junction21_test2_iE1.txt};
		\addlegendentry{Burgers}
		\end{axis}
		
		\begin{axis}[ylabel = $u^2$,
		y=2cm,
		width=8cm,
		xmin = -0.0, xmax = 1.0,
		tick label style={/pgf/number format/fixed, 
			/pgf/number format/precision=5},
		scaled x ticks=false,
		ymin = -1.1, ymax = 0.1,
		name=edge2,
		at=(edge1.below south west), anchor=above north west,
		]
		\addplot[color = blue,thick] file{Data/rho_kinetic_junction21_test2_iE2.txt};
		\addplot[color = red,thick] file{Data/rho_Burgers_junction21_test2_iE2.txt};
		\end{axis}
		
		\begin{axis}[ylabel = $u^3$,xlabel =  $x$,
		y=2cm,
		width=8cm,
		xmin = -0.0, xmax = 1.0,
		tick label style={/pgf/number format/fixed, 
			/pgf/number format/precision=5},
		scaled x ticks=false,
		ymin = -1.1, ymax = 0.1,
		at=(edge2.below south west), anchor=above north west,
		]
		\addplot[color = blue,thick] file{Data/rho_kinetic_junction21_test2_iE3.txt};
		\addplot[color = red,thick] file{Data/rho_Burgers_junction21_test2_iE3.txt};
		\end{axis}
		\end{tikzpicture}
	}
	\caption{Left: RP$1$-$1$-$1$ with $u^1=-1$,  $u^2=-0.75$,  $u^3=0.5$, Right: RP$1$-$1$-$2$ with $u^1=-1$,   $u^2=-0.75$,  $u^3=-0.5$.}
	\label{fig:Network_J21_111_112}
\end{figure}
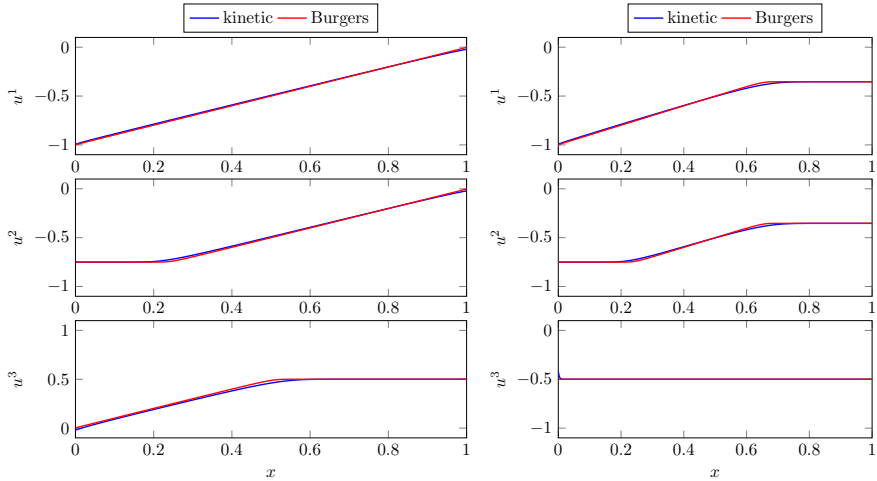
	
Figure \ref{fig:Network_J21_111_112} shows Riemann Problems to the cases RP$1$-$1$-$1$ and RP$1$-$1$-$2$.
On the left the first two rarefaction waves move to the left, the third one to the right.
On the right hand side the flow from the third edge is distributed equally to the first two edges.
	
	\begin{figure}
		\externaltikz{junction21_121_122}{
			\begin{tikzpicture}[scale=0.65]
			\begin{axis}[ylabel = $u^1$,
			y=2cm,
			x=8cm,
			xmin = -0.0, xmax = 1.0,
			ymin = 0.7, ymax = 1.5,
			legend style = {at={(0.5,1)},xshift=0.2cm,yshift=0.1cm,anchor=south},
			legend columns= 3,
			name=edge1,
			]
			\addplot[color = blue,thick] file{Data/rho_kinetic_junction21_test3_iE1.txt};
			\addlegendentry{kinetic}
			\addplot[color = red,thick] file{Data/rho_Burgers_junction21_test3_iE1.txt};
			\addlegendentry{Burgers}
			\end{axis}
			
			\begin{axis}[ylabel = $u^2$,
			y=2cm,
			x=8cm,
			xmin = -0.0, xmax = 1.0,
			ymin = -0.9, ymax = 1.1,
			name=edge2,
			at=(edge1.below south west), anchor=above north west,
			]
			\addplot[color = blue,thick] file{Data/rho_kinetic_junction21_test3_iE2.txt};
			\addplot[color = red,thick] file{Data/rho_Burgers_junction21_test3_iE2.txt};
			\end{axis}
			
			\begin{axis}[ylabel = $u^3$,xlabel =  $x$,
			y=2cm,
			x=8cm,
			xmin = -0.0, xmax = 1.0,
			ymin = 0.3, ymax = 1.1,
			at=(edge2.below south west), anchor=above north west,
			]
			\addplot[color = blue,thick] file{Data/rho_kinetic_junction21_test3_iE3.txt};
			\addplot[color = red,thick] file{Data/rho_Burgers_junction21_test3_iE3.txt};
			\end{axis}
			\end{tikzpicture}
			\begin{tikzpicture}[scale=0.65]
			\begin{axis}[ylabel = $u^1$,
			y=2cm,
			y=2cm,
			width=8cm,
			xmin = -0.0, xmax = 1.0,
			scaled x ticks=false,
			ymin = -1.1, ymax = 0.1,
			legend style = {at={(0.5,1)},xshift=0.2cm,yshift=0.1cm,anchor=south},
			legend columns= 3,
			name=edge1,
			]
			\addplot[color = blue,thick] file{Data/rho_kinetic_junction21_test4_iE1.txt};
			\addlegendentry{kinetic}
			\addplot[color = red,thick] file{Data/rho_Burgers_junction21_test4_iE1.txt};
			\addlegendentry{Burgers}
			\end{axis}
			
			\begin{axis}[ylabel = $u^2$,
			y=2cm,
			width=8cm,
			xmin = -0.0, xmax = 1.0,
			ytick distance=0.3,
			scaled x ticks=false,
			ymin = -0.6, ymax = 0.6,
			name=edge2,
			at=(edge1.below south west), anchor=above north west,
			]
			\addplot[color = blue,thick] file{Data/rho_kinetic_junction21_test4_iE2.txt};
			\addplot[color = red,thick] file{Data/rho_Burgers_junction21_test4_iE2.txt};
			\end{axis}
			
			\begin{axis}[ylabel = $u^3$,xlabel =  $x$,
			y=2cm,
			width=8cm,
			xmin = -0.0, xmax = 1.0,
			tick label style={/pgf/number format/fixed, 
				/pgf/number format/precision=5},
			scaled x ticks=false,
			ymin = -1.1, ymax = 0.1,
			at=(edge2.below south west), anchor=above north west,
			]
			\addplot[color = blue,thick] file{Data/rho_kinetic_junction21_test4_iE3.txt};
			\addplot[color = red,thick] file{Data/rho_Burgers_junction21_test4_iE3.txt};
			\end{axis}
			\end{tikzpicture}
		}
		\caption{Left: RP$2$-$1$-$1$ with $u^1=1$,  $u^2=-0.75$,  $u^3=0.5$, Right: RP$1$-$2$-$2$ with $u^1=-1$,   $u^2=0.5$,  $u^3=-0.6$.}
		\label{fig:Network_J21_121_122}
	\end{figure}
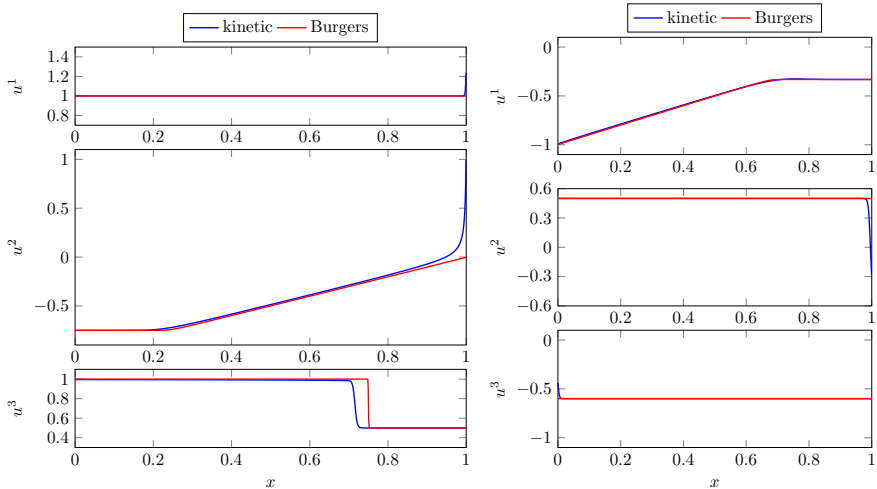
	On the left side of Figure \ref{fig:Network_J21_121_122} the flow enters from the first edge and is directed completely into the third one. 
	In edge $2$ the flow at the junction is zero, which leads to a rarefaction wave and a strong layer. 
	Note that this transition connects states with negative and positive sign.
	In the third edge the wave in the kinetic model is behind the macroscopic shock due to the temporal layer from the initial conditions. 
	This gap closes when decreasing $\epsilon$ and the grid spacing. 
	On the right hand side the flow from edge $2$ and $3$ enters into edge $1$, which corresponds to the second subcase of RP$1$-$2$-$2$.
	
	\begin{figure}
		\externaltikz{junction21_221_222}{
			\begin{tikzpicture}[scale=0.65]
			\begin{axis}[ylabel = $u^1$,
			y=2cm,
			x=8cm,
			xmin = -0.0, xmax = 1.0,
			ymin = -0.1, ymax = 1.1,
			legend style = {at={(0.5,1)},xshift=0.2cm,yshift=0.1cm,anchor=south},
			legend columns= 3,
			name=edge1,
			]
			\addplot[color = blue,thick] file{Data/rho_kinetic_junction21_test5_iE1.txt};
			\addlegendentry{kinetic}
			\addplot[color = red,thick] file{Data/rho_Burgers_junction21_test5_iE1.txt};
			\addlegendentry{Burgers}
			\end{axis}
			
			\begin{axis}[ylabel = $u^2$,
			y=2cm,
			x=8cm,
			xmin = -0.0, xmax = 1.0,
			ymin = -.1, ymax = 1.1,
			name=edge2,
			at=(edge1.below south west), anchor=above north west,
			]
			\addplot[color = blue,thick] file{Data/rho_kinetic_junction21_test5_iE2.txt};
			\addplot[color = red,thick] file{Data/rho_Burgers_junction21_test5_iE2.txt};
			\end{axis}
			
			\begin{axis}[ylabel = $u^3$,xlabel =  $x$,
			y=2cm,
			x=8cm,
			xmin = -0.0, xmax = 1.0,
			ymin = -.1, ymax = 1.1,
			at=(edge2.below south west), anchor=above north west,
			]
			\addplot[color = blue,thick] file{Data/rho_kinetic_junction21_test5_iE3.txt};
			\addplot[color = red,thick] file{Data/rho_Burgers_junction21_test5_iE3.txt};
			\end{axis}
			\end{tikzpicture}
			\begin{tikzpicture}[scale=0.65]
			\begin{axis}[ylabel = $u^1$,
			y=2cm,
			y=2cm,
			width=8cm,
			xmin = -0.0, xmax = 1.0,
			scaled x ticks=false,
			ymin = -0.0, ymax = 1.1,
			legend style = {at={(0.5,1)},xshift=0.2cm,yshift=0.1cm,anchor=south},
			legend columns= 3,
			name=edge1,
			]
			\addplot[color = blue,thick] file{Data/rho_kinetic_junction21_test6_iE1.txt};
			\addlegendentry{kinetic}
			\addplot[color = red,thick] file{Data/rho_Burgers_junction21_test6_iE1.txt};
			\addlegendentry{Burgers}
			\end{axis}
			
			\begin{axis}[ylabel = $u^2$,
			y=2cm,
			width=8cm,
			xmin = -0.0, xmax = 1.0,
			tick label style={/pgf/number format/fixed, 
				/pgf/number format/precision=5},
			scaled x ticks=false,
			ymin = -0.0, ymax = 1.1,
			name=edge2,
			at=(edge1.below south west), anchor=above north west,
			]
			\addplot[color = blue,thick] file{Data/rho_kinetic_junction21_test6_iE2.txt};
			\addplot[color = red,thick] file{Data/rho_Burgers_junction21_test6_iE2.txt};
			\end{axis}
			
			\begin{axis}[ylabel = $u^3$,xlabel =  $x$,
			y=2cm,
			width=8cm,
			xmin = -0.0, xmax = 1.0,
			tick label style={/pgf/number format/fixed, 
				/pgf/number format/precision=5},
			scaled x ticks=false,
			ymin = -.6, ymax = 0.8,
			at=(edge2.below south west), anchor=above north west,
			]
			\addplot[color = blue,thick] file{Data/rho_kinetic_junction21_test6_iE3.txt};
			\addplot[color = red,thick] file{Data/rho_Burgers_junction21_test6_iE3.txt};
			\end{axis}
			\end{tikzpicture}
		}
		\caption{Left: RP$2$-$2$-$1$ with $u^1=0.5$,  $u^2=0.4$,  $u^3=0.75$, Right: RP$2$-$2$-$2$ with $u^1=0.5$,   $u^2=0.4$,  $u^3=-0.3$.}
		\label{fig:Network_J21_221_222}
	\end{figure}
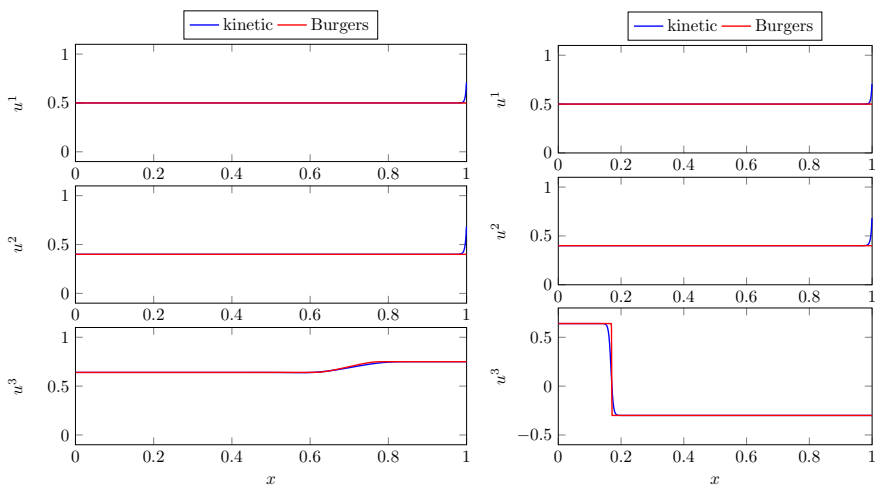
	The final two tests are shown in Figure \ref{fig:Network_J21_221_222}.
	The flow from the first two edges merges into the third one. 
	On the right hand side  this leads to a strong transsonic shock wave. 	
	In this test again a linearization would fail, as the characteristics in one edge have to be flipped to guarantee the conservation of mass. 
	
	Also for the $2$-$1$-Junction all tests show a good agreement of the kinetic and the macroscopic solution. 
	The states at the junction are correctly identified by the results of section \ref{sec:21-node}.

\section{Conclusion and Outlook}
Starting from a kinetic network model with prescribed coupling conditions at the nodes we have derived coupling conditions for a nonlinear scalar conservation law.
This is achieved by identifying the asymptotic states and their domains of attraction in the kinetic layer problems at the nodes and combining the layer solutions  with  half Riemann-problems for the macroscopic variables. 
The resulting coupling conditions conserve the total mass and can handle arbitrary flow scenarios.
The results are compared numerically to the full solution of the kinetic problem on the network.

The presented approach can be extended to more general  flux functions $F=F(u)$ or other choices of kinetic coupling conditions.  
The case of nonlinear hyperbolic systems is subject to ongoing investigations with a similar procedure.
	
\section*{ Acknowledgment}
This work has been supported by Deutsche Forschungsgemeinschaft (DFG) with the grant BO 4768/1.




\begin{thebibliography}{10}
	\bibliographystyle{siam}
	
	\bibitem{AS15}{B. Andreianov, K. Sbihi,
		{\em Well-posedness of general boundary-value problems for scalar conservation laws},
		{Trans. Amer. Math. Soc.}, 367(6), 3763--3806, 2015.
	}	
	\bibitem{ACD17}{B.P. Andreianov, G.M. Coclite, C. Donadello,
		{\em Well-posedness for vanishing viscosity solutions of scalar conservation laws on a network},
		Discrete Contin. Dyn. Syst., 37(11), 5913--5942, 2017
	}
	\bibitem{AM04} 
	D. Aregba-Driollet,V. Milisic,
	{\em Kinetic approximation of a boundary value problem for conservation laws},
	Numer. Math.  97, 595–633, 2004
	\bibitem{BHK06a}
	M. Banda, M. Herty, A. Klar, {\em Gas flow in pipeline networks}, NHM 1(1), 41-56, 2006	
	\bibitem{BRN}{C. Bardos, A.Y. le Roux,J.-C. N\'ed\'elec, 
		{\em First order quasilinear equations with boundary conditions},
		Comm. Partial Differential Equations, 4(9), 1017--1034, 1979.
		}		
	\bibitem{BSS84} C. Bardos, R. Santos, and R Sentis, {\em Diffusion approximation and computation of the critical size}, Trans. Amer. Math. Soc. 284, 2, 617-649, 1984
	\bibitem{BLP79} A. Bensoussan, J.L. Lions, and G.C. Papanicolaou, {\em Boundary-layers and homogenization of transport processes}, J. Publ. RIMS Kyoto Univ. 15, 53-157, 1979
	\bibitem{BK17} R. Borsche, A. Klar {\em  Kinetic derived coupling conditions for macroscopic equations},arxiv 2017
	\bibitem{BKP16} R. Borsche, A. Klar, T.N.H. Pham, {\em  Kinetic and related macroscopic models for chemotaxis on networks}, M3AS, 26, No. 6, 1219-1242, 2016
	\bibitem{BNR14}{ G. Bretti, R. Natalini, M. Ribot,
		{\em A hyperbolic model of chemotaxis on a network: a numerical study},
		ESAIM Math. Model. Numer. Anal., 48(1) ,231--258, 2014.
	}		
	\bibitem{BFK07}{R. B\"urger, H. Frid, K.H. Karlsen
		{\em On the well-posedness of entropy solutions to conservation laws with a zero-flux boundary condition}, J. Math. Anal. Appl., 326(1), 108--120, 2007.
	}	
	\bibitem{C69} C. Cercignani, {\em A Variational Principle for Boundary Value Problems},
	J. of Stat. Phys., Vol 1, No. 2, 1969
	\bibitem{C88} C. Cercignani, {\em The Boltzmann Equation and its Applications},
	Springer, 1988
	\bibitem{CLT14} I.K. Chen, T.P. Liu, and S. Takata, {\em Boundary singularity for thermal transpiration problem of the linearized Boltzmann equation}, Arch. Ration. Mech. Anal. 212, 2, 575–595, 2014. 
	\bibitem{CGP05}
	{G.~M. Coclite, M.~Garavello, and B.~Piccoli}, {\em Traffic flow on a road
		network}, SIAM J. Math. Anal., 36 , 1862--1886 2005.	
	\bibitem{CG08}{
		R.M. Colombo, M. Garavello,
		{\em On the Cauchy problem for the $p$-system at a junction},
		SIAM J. Math. Anal., 39, 1456--1471 2008.
	}	
	\bibitem{CM08}{
		R.M. Colombo, Rinaldo,C. Mauri,
		{\em Euler system for compressible fluids at a junction},
		J. Hyperbolic Differ. Equ., 5(3), 547--568, 2008.
	}	
	\bibitem{Cor17}
	A. Corli, L. di Ruvo, L. Malaguti, M. D. Rosini,
	Traveling waves for degenerate diffusive equations on networks, NHM 12,3, 339 - 370, 2017
	\bibitem{Coron} F. Coron, {\em Computation of the Asymptotic States for Linear Halfspace
		Problems}, TTSP 19(2),  89, 1990
	\bibitem{CGS} F. Coron, F. Golse, C. Sulem, {\em A Classification of Well-posed Kinetic
		Layer Problems}, CPAM, Vol. 41,  409, 1988
	\bibitem{Zuazua2} R. Dager, E. Zuazua, {\em Controllability of tree-shaped networks of vibrating strings}, C. R. Acad. Sci. Paris, 332, I, 1087–1092, 2001
	\bibitem{G08} F. Golse, {\em Analysis of the boundary layer equation in the kinetic theory of gases}, Bull. Inst. Math. Acad. Sin. 3, 1, 211-242, 2008
	\bibitem{GK95} F. Golse, A. Klar, {\em Numerical Method for Computing Asymptotic States and Outgoing Distributions for a Kinetic Linear Half Space Problem}, J. Stat. Phys. 80 (5-6), 1033-1061, 1995
	\bibitem{GMP} W. Greenberg, C. van der Mee, V. Protopopescu, {\em Boundary Value
		Problems in Abstract Kinetic Theory}, Birkh\"auser, 1987	
	\bibitem{GPBook} M. Garavello, B. Piccoli, {\em Traffic flow on networks}, AIMS, 2006
	\bibitem{HM09}
	M.~Herty and S.~Moutari.
	\newblock A macro-kinetic hybrid model for traffic flow on road networks.
	\newblock {\em Comput. Methods Appl. Math.}, 9, 3,238--252, 2009.	
	\bibitem{LS02}{
		G. Leugering, Guenter, E.J.P.G. Schmidt,
		{\em On the modelling and stabilization of flows in networks of open canals},
		SIAM J. Control Optim., 41(1), 164--180 2002.
	}	
	\bibitem{Liu} T.P. Liu, {\em Hyperbolic conservation laws with relaxation}, Commun. Math. Phys. 108, 153-175 (1987) 
	\bibitem{LLS162} Q. Li, J. Lu, and W. Sun, {\em Half-space kinetic equations with general boundary conditions}, Math. Comp. 2016
	\bibitem{LLS16} Q. Li,J.Lu,W. Sun, {\em A convergent method for linear half-space kinetic equations}, arxiv, 2015
	\bibitem{LY01}
	H. Liu and W.-A. Yong, {\em Time-asymptotic stability of boundary-layers for a hyperbolic relaxation system}, Comm. Partial Differential Equations, 26(7-8), 1323–1343, 2001.
	\bibitem{LX96}
	J.-G. Liu, Z. Xin, {\em Boundary-layer behavior in the fluid-dynamic limit for a
	nonlinear model Boltzmann Equation}, Arch. Rational Mech. Anal. 135,
	61-105, 1996.
	\bibitem{RT01}
	R. Natalini and A. Terracina, {\em Convergence of a relaxation approximation to a boundary value problem for conservation laws}, Comm. Partial Differential Equations, 26(7-8), 1235–1252, 2001.
	\bibitem{N96}
	S. Nishibata, {\em The initial boundary value problems for hyperbolic conservation laws
	with relaxation}, J. Diff. Eqns. 130, 100-126, 1996.
	\bibitem{N97}
	S. Nishibata, S.-H. Yu, {\em The asymptotic behavior of the hyperbolic conservation
	laws with relaxation on the quarter-plane}, Siam J. Math. Anal. 28 , 304-321, 1997.
	\bibitem{O96}{F. Otto, {\em Initial-boundary value problem for a scalar conservation law},
		C. R. Acad. Sci. Paris S\'er. I Math., 322(8), 729--734, 1996.
	}
	\bibitem{SO} Y. Sone, Y. Onishi, {\em Kinetic Theory of Evaporation and Condensation,
		Hydrodynamic Equation and Slip Boundary Condition}, J. Phys. Soc. of Japan, Vol.
	44, No. 6,  1981, 1978
	\bibitem{UTY03} S. Ukai, T. Yang, and S.-H. Yu, {\em Nonlinear boundary layers of the Boltzmann equation. I. Existence}, Comm. Math. Phys. 236,  3, 373-393, 2003
	\bibitem{ToroBook}E.F. Toro {\em Riemann solvers and numerical methods for fluid dynamics}, Springer, 2009
	\bibitem{WX99}
	W.-C. Wang, Z. Xin,
	{\em Asymptotic limit of initial boundary value problems for conservation laws with relaxational extensions},
	Communications on Pure and Applied Mathematics, 51,5  505–535, 1998
	\bibitem{X04}
	W.-Q. Xu, 
	Boundary conditions and boundary layers for a multi-dimensional relaxation model,
	Journal of Differential Equations 197,  1, 10, 85-117, 2004
	\bibitem{WY99}
	Wen-An Yong, {\em Boundary conditions for hyperbolic systems with stiff relaxation}, Indiana University Mathematics Journal 48, 1, 115-137,
	1999
	\bibitem{Zuazua1} J. Valein, E. Zuazua, {\em Stabilization of the Wave Equation on 1-d Networks}, SIAM Journal on Control and Optimization, 48, 4, 2771-2797, 2009 
	
	
\end{thebibliography}
\end{document}